\documentclass[12pt]{article}
\usepackage{mathrsfs}
\usepackage{amsmath}
\usepackage{amssymb}
\usepackage{amsfonts}
\usepackage{graphicx}
\usepackage{wasysym}
\usepackage{verbatim}
\usepackage[amsmath,thmmarks]{ntheorem}

\def\F{\mathscr{F} }

\def\R{\mathbb{R}}
\def\Z{\mathbb{Z}}
\def\T{\mathbb{T}}

\def\ba{\begin{array}}
\def\ea{\end{array}}
\def\be{\begin{enumerate}}
\def\ee{\end{enumerate}}
\def\bi{\begin{itemize}}
\def\ei{\end{itemize}}
\def\bc{\begin{cases}}
\def\ec{\end{cases}}
\def\beq{\begin{equation}}
\def\eeq{\end{equation}}
\def\beqs{\begin{equation*}}
\def\eeqs{\end{equation*}}
\def\beqa{\begin{eqnarray}}
\def\eeqa{\end{eqnarray}}
\def\beqas{\begin{eqnarray*}}
\def\eeqas{\end{eqnarray*}}
\def\bmul{\begin{multline}}
\def\emul{\end{multline}}
\def\bmuls{\begin{multline*}}
\def\emuls{\end{multline*}}
\def\bg{\begin{gather}}
\def\eg{\end{gather}}
\def\bgs{\begin{gather*}}
\def\egs{\end{gather*}}

\theorembodyfont{\normalfont\sl}
\newtheorem{thm}{Theorem}[section]
\newtheorem{cor}[thm]{Corollary}
\newtheorem{lem}[thm]{Lemma}
\newtheorem{prop}[thm]{Proposition}

\numberwithin{equation}{section}

\begin{document}

\leftskip 0 true cm \rightskip 0 true cm

\newpage

\setcounter{page}{1}
 \baselineskip=24pt

\title{
\baselineskip=24pt
 \bf The Cauchy Problem of the\\[.5ex]
 \bigskip
 \bf Schr\"{o}dinger-Korteweg-de Vries System\footnote{This work is
 supported by National Natural Science Foundation of China
 under grant numbers 10471047 and 10771074.}}

\bigskip

\author{
\baselineskip=24pt
 \normalsize Yifei WU\footnote{Email: yerfmath@yahoo.cn}\\[1ex]
 \normalsize Department of Mathematics,
             South China University of Technology, \\[1ex]
 \normalsize Guangzhou, Guangdong 510640, P. R. China
 }

\bigskip

\date{}

\bigskip

\bigskip

\leftskip 1cm \rightskip 1cm

\maketitle

\noindent
 {\small{\bf Abstract}\quad
We study the Cauchy problem of the Schr\"{o}dinger-Korteweg-de
Vries system. First, we establish the local well-posedness
results, which improve the results of Corcho, Linares (2007). Moreover,
we obtain some ill-posedness results, which show
that they are sharp in some well-posedness thresholds.
Particularly, we obtain the local well-posedness for the initial data in
$H^{-\frac{3}{16}+}(\R)\times H^{-\frac{3}{4}+}(\R)$ in the resonant case,
it is almost the optimal except the endpoint.
At last we establish the global well-posedness results
in $H^s(\R)\times H^s(\R)$ when $s>\dfrac{1}{2}$
no matter in the resonant case or in the non-resonant case, which improve the
results of Pecher (2005). }
\bigskip

\noindent
 {\small{\bf Keywords:}\quad Schr\"{o}dinger-Korteweg-de Vries system,
local well-posedness,
ill-posedness, global well-posednsss, Bourgain space, $I$-method}
\bigskip

\noindent
 {\small{\bf  MR(2000) Subject Classification:}\quad 35Q53, 35Q55}

\bigskip

\bigskip

\leftskip 0cm \rightskip 0cm

\normalsize

\baselineskip=24pt

\section{Introduction}
The Cauchy problem of the Schr\"{o}dinger-Korteweg-de Vries equations
\renewcommand{\arraystretch}{2}
\begin{equation}\label{NLS-KdV}
\left\{
\begin{array}{l}
  i\partial_{t}u+\partial_{x}^{2}u=\alpha uv+ \beta |u|^2u,
     \qquad x,\,t\in \R,\\
  \partial_{t}v+\partial_{x}^{3}v+\dfrac{1}{2}\partial_xv^2
  =\gamma \partial_x(|u|^2),\\
  u(x,0)=u_0(x)\in H^s(\R), v(x,0)=v_0(x)\in H^l(\R),
\end{array}
\right.
\end{equation}
where $\alpha,\beta,\gamma\in \R$. The system governs the
interactions between the short-wave and the long-wave, which appears
in several fields of physics and fluid dynamics. The case $\beta=0$
describes the resonant interactions, while the case $\beta\neq0$
describes the non-resonant interactions. See \cite{FO},
 \cite{HIMN}, \cite{KKS}, \cite{SY} for the applications.

The Cauchy problem for the system (\ref{NLS-KdV}) was considered by several
authors. The local well-posedness was studied in \cite{BOP}, \cite{BOP2},
\cite{CL}, $etc$., where the last paper \cite{CL} obtained the local well-posedness
for $(u_0,v_0)\in H^s(\R)\times H^l(\R)$ when
$s\geq 0$, $l>-\dfrac{3}{4}$, and
\begin{itemize}
        \item $s-1\leq l \leq 2s-1/2,$ if $s\leq 1/2$;
        \item $s-1\leq l < s+1/2,$ if $s> 1/2$,
\end{itemize}
by the Bourgain argument (see \cite{B}, \cite{KPV1} for instances).
Moreover, when $\alpha\gamma>0$, Tsutsumi \cite{T} proved by some
conservation laws that for $(u_0,v_0)\in H^{s+\frac{1}{2}}(\R)
\times H^{s}(\R)$ with $s\in \Z^+$, (\ref{NLS-KdV}) was global
well-posedness; Guo and Miao \cite{GM} showed that the system in the
resonant case was globally well-posed for $(u_0,v_0)\in H^{s}(\R)
\times H^{s}(\R)$ with $s\in \Z^+$; In \cite{P}, the author improved
the results and obtained the global well-posedness for $(u_0,v_0)\in
H^{s}(\R)\times H^{s}(\R)$ when $s<1$ and
\begin{itemize}
        \item $s>3/5$ in the case of $\beta=0$;
        \item $s>2/3$ in the case of $\beta\neq0$,
\end{itemize}
by using the $I-$method of Colliander, Keel, Staffilani,  Takaoka
and Tao (see \cite{CKSTT}, \cite{CKSTT2} for examples).

In \cite{CL}, the authors have obtained the local well-posedness
for the initial data belongs to $L^2(\R)\times H^{-\frac{3}{4}+}(\R)$,
but it seems not the natural one as the best result
(except the endpoint) and also exists some room
if the system (\ref{NLS-KdV}) has no
power type nonlinearity (that is, $\beta=0$). As what studied
in the first part of this paper, we establish the local well-posedness
results at a relatively wide region of the indices $(s,l)$, compared to
the results in \cite{CL}.
Especially, we obtain the local well-posedness
for $(u_0,v_0)\in H^{-\frac{3}{16}+}(\R)\times H^{-\frac{3}{4}+}(\R)$
in the case of $\beta=0$.
The second aim here is to establish some ill-posedness results,
by the breakage of continuity in Picard iterative scheme
(see \cite{CCT1}, \cite{MR}, \cite{BT}, $etc$.), which show that
some thresholds are sharp except the boundary
in the well-posedness region. By these results, we will see that
the index $(-\dfrac{3}{16}+,-\dfrac{3}{4}+)$ is almost the best in the
resonant case. Further, the third aim in this article is to obtain the
global well-posedness by the $I$-method. As we know, the thresholds of
the global well-posedness in the Sobolev spaces, studied by the $I$-method,
are decided by two ingredients: the almost conserved
quantities and the lifetime in the local theory, particularly
if the solutions of the equations lack of the scale invariance.
Our motivation here is to lengthen the lifetime of the local existence.
Our argument is establishing some special type multilinear estimates,
as an available technique we obtain a uniformly control of
the Bourgain $X_{0,b}-$norm of $u$ by the $L^2$-mass conservation.
But unfortunately, these special type estimates break the framework of
the fixed point theorem, we finally use the iterate technical to overcome
this problem. However, we believe that our global results may not be the best
and might be improved especially
by some more sophisticated estimates on the almost conserved quantities.

\noindent{\bf Some basic notations}. We use $A\lesssim B$ or $B\gtrsim A$
to denote the statement that $A\leq CB$ for some large constant $C$
which may vary from line to line.  We use $A\ll
B$ to denote the statement $A\leq C^{-1}B$, and use $A\sim B$ to
mean $A\lesssim B\lesssim A$. The notation $a+$ denotes $a+\epsilon$
for any small $\epsilon$, and $a-$ for $a-\epsilon$.
$\langle\cdot\rangle=(1+|\cdot|^2)^{\frac{1}{2}}$  and
$D_x=(-\partial^2_x)^{\frac{1}{2}}$. We use $\|f\|_{L^p_xL^q_t}$
to denote the mixed norm
$\Big(\displaystyle\int\|f(x,\cdot)\|_{L^q}^p\
dx\Big)^{\frac{1}{p}}$. Moreover, we denote
$\hat{u}$ to be the spatial or
spacetime Fourier transform of
$u$, and use $\check{f}$ or $\F^{-1}$ (such as
$\F_{\xi}^{-1}, \F_{\tau}^{-1}, \F_{\xi\tau}^{-1}$ $ect$.) to
denote the inverse Fourier transform of $f$ (on the corresponding variables).

Now we introduce some definitions before presenting our main
results.
We use $U_\phi(t)= \exp(-it \phi(-i\partial_x))$ to denote the
unitary group generated by the linear equation
$$
iu_t-\phi(-i\partial_x)u=0,
$$
and define  the Bourgain spaces $X_{s,b}(\phi)$ to be the closure of
the Schwartz class under the norms
$$
\|f\|_{X_{s,b}(\phi)}\equiv\left(\displaystyle\int\!\!\!\!\int
\langle\xi\rangle^{2s}\langle\tau+\phi(\xi)\rangle^{2b}|\hat{f}(\xi,\tau)|^2
\,d\xi d\tau\right)^{\frac{1}{2}},
$$
for $s,b\in \R $. We write $X_{s,b}^{\pm}, Y_{s,b}$
to be $X_{s,b}(\phi)$ when $\phi=\pm\xi^2$,
$-\xi^3$, which is corresponding to the Schr\"{o}dinger and KdV respectively,
and we write $X_{s,b}\equiv X_{s,b}^{+}$ in default.
For an interval $\Omega$, we define $X_{s,b}^{\Omega}(\phi)$
to be the restriction of $X_{s,b}(\phi)$ on $\R\times\Omega$
with the norms
$$
\|f\|_{X_{s,b}^{\Omega}(\phi)}
=\inf\{\|F\|_{X_{s,b}(\phi)}:F|_{t\in\Omega}=f|_{t\in\Omega}\}.
$$
When $\Omega=[-\delta,\delta]$,
we write $X_{s,b}^\Omega(\phi)$ as $X_{s,b}^\delta(\phi)$
($X_{s,b}^\Omega$ as $X_{s,b}^\delta$,
$Y_{s,b}^\Omega$ as $Y_{s,b}^\delta$).

Let $s<0$ and $N\gg 1$ be fixed,  the Fourier
multiplier operator $I_{N,s}$ is defined as
\beq
\widehat{I_{N,s}u}(\xi)=m_{N,s}(\xi)\hat{u}(\xi),\label{I}
\eeq
where the multiplier $m_{N,s}(\xi)$ is a smooth, monotone function
satisfying $0<m_{N,s}(\xi)\leq 1$ and
 \beq m_{N,s}(\xi)=\Biggl\{
\begin{array}{ll}
1,&|\xi|\leq N,\\
N^{1-s}|\xi|^{s-1},&|\xi|>2N.\label{m}
\end{array}
\eeq
Sometimes we denote $I_{N,s}$ and $m_{N,s}$ as $I$ and $m$
respectively for short if there is no confusion.
It is obvious that the operator $I_{N,s}$ maps $H^s(\R)$ into
$H^1(\R)$ with equivalent norm such that
\beq
    \|f\|_{H^s}
\lesssim
    \|I_{N,s}f\|_{H^1}
\lesssim
    N^{1-s}\|f\|_{H^s}.\label{II}
\eeq
Moreover, $I_{N,s}$ can be extended to a map (still denoted by
$I_{N,s}$) from $X_{s,b}$ to $X_{1,b}$ which satisfies
$$
    \|f\|_{X_{s,b}}\lesssim\|I_{N,s}f\|_{X_{1,b}}
\lesssim N^{1-s}\|f\|_{X_{s,b}}
$$
for any $s<1,b \in \R$.

Our main results in this paper are given as follows.

\begin{thm}
The Cauchy problem of the system (\ref{NLS-KdV}) is locally
well-posed on some time interval $[-\delta, \delta]$ for the initial
data $(u_0,v_0)\in H^s(\R)\times H^l(\R)$ when
\begin{itemize}
        \item $l>-3/4, l<4s, s-2\leq l < s+1,$ if $\beta=0$;
        \item $l>-3/4, s\geq 0, l<4s, s-2\leq l < s+1,$ if $\beta\neq0$.
\end{itemize}
The solutions satisfy
$$
(u,v)\in C_t^0([-\delta, \delta]; H^s(\R)\times H^l(\R)).
$$
\end{thm}

{\noindent\it Remark.} The best result obtained in Theorem 1.1 is local
well-posedness in $L^2(\R)\times H^{-\frac{3}{4}+}(\R)$ in the
non-resonant case ($\beta\neq0$), which has been contained in \cite{CL}.
However, in the resonant case ($\beta=0$),
the best result we obtained is local
well-posedness in $H^{-\frac{3}{16}+}(\R)\times H^{-\frac{3}{4}+}(\R)$,
which improves the one in \cite{CL}.
As we see, it follows from the assumptions of $l>-3/4$ and $l<4s$. The next
result tells that it is almost the best in the
following sense except the endpoint case.
\begin{thm}
Let $l>4s$, and let $\delta$ and
the solution map $(u_0,v_0)\mapsto (u,v)$ is defined in
Theorem 1.1 from $H^{s_0}(\R)\times H^{l_0}(\R)$ to
$C_t^0([0,\delta];H^{s_0}(\R)\times H^{l_0}(\R))$
for some well-posed index $(s_0,l_0)$.
Then the map is not $C^2$-differentiable at zero from
$H^{s}(\R)\times H^{l}(\R)$ to
$C_t^0([0,\delta];H^{s}(\R)\times H^{l}(\R))$.
\end{thm}

If we take the initial data $(0,v_0)$ with $v_0\in H^l(\R)$
in the system (\ref{NLS-KdV}),
then by the uniqueness, $u\equiv0$ and the  system (\ref{NLS-KdV})
can be deduced to the single KdV equation which is well known
to ill-posedness when $l<-\dfrac{3}{4}$ (see \cite{CCT}, \cite{KPV2}).
Thus we say that the system (\ref{NLS-KdV}) is ill-posedness
in $H^s(\R)\times H^l(\R)$ for $s\in \R$ and $l<-\dfrac{3}{4}$.
To sum up, we draw the following figure for the corresponding regions.

\centerline{\includegraphics[scale=1, trim=0 15cm 0 1cm]
{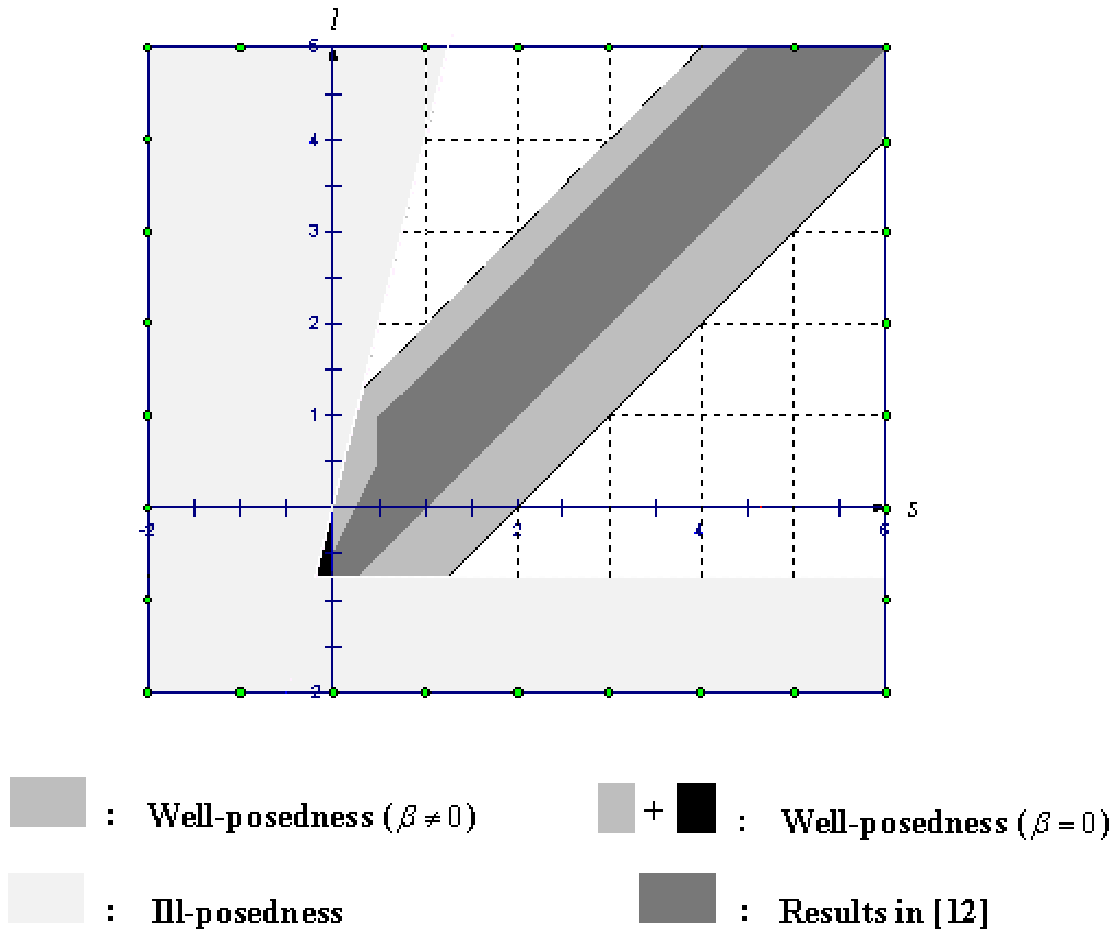}} \centerline{\bf {Figure 1: Well-posedness and
ill-posedness for the NLS-KDV system}}

\begin{thm}
Let $\alpha\gamma>0$, then the Cauchy problem of
the system (\ref{NLS-KdV}) is globally well-posed for the
initial data $(u_0,v_0)\in H^s(\R)\times H^s(\R)$ when
$1>s>\dfrac{1}{2}$, no matter in the resonant case or in the non-resonant case.
\end{thm}

The rest of this article is organized as follows.
In Section 2, we derive some preliminary estimates.
In Section 3, we establish some multilinear estimates and prove Theorem 1.1.
In Section 4, we establish some multilinear estimates of special types,
give a variant local result, and prove Theorem 1.3.
In Section 5, we prove Theorem 1.2.
In Section 6, as an appendix, we prove some auxiliary lemmas about the spaces
$X_{s,b}^\Omega(\phi)$.

 \vspace{0.3cm}
\section{Some Preliminary Estimates}

First, we present two Stricharz estimates in the Bourgain
spaces.
\begin{lem}
(1) For any $u\in X_{0,\theta+}^{\pm}$,
$\theta\geq  \dfrac{3}{2}\left(\dfrac{1}{2}-\dfrac{1}{p}\right)$, $p\in [2,6]$,
we have
\begin{equation}
\|u\|_{L^p_{xt}}
    \lesssim
    \|u\|_{X_{0,\theta+}^{\pm}}.
    \label{XE1}
\end{equation}
(2) For any $v\in Y_{0,\rho+}$,
$\rho\geq  \dfrac{4}{3}\left(\dfrac{1}{2}-\dfrac{1}{q}\right)$.
$q\in [2,8]$, we have
\begin{equation}
\|v\|_{L^q_{xt}}
    \lesssim
    \|v\|_{Y_{0,\rho+}},
    \label{XE2}
\end{equation}
\end{lem}
\begin{proof}
They are easily obtained by the interpolation between
the following well-known inequalities
$$
\|u\|_{L^6_{xt}}
    \lesssim
    \|u\|_{X_{0,\frac{1}{2}+}};\quad
\|v\|_{L^8_{xt}}
    \lesssim
    \|v\|_{Y_{0,\frac{1}{2}+}}
$$
and the equalities
$$
\|u\|_{L^2_{xt}}
    =
    \|u\|_{X_{0,0}};\quad
\|v\|_{L^2_{xt}}
    =
    \|v\|_{Y_{0,0}}.
$$
\hfill$\Box$
\end{proof}

Now, we introduce some multiplier operators and the estimates on them
(similar results are appeared in \cite{CKS2} and \cite{G}).
They are the important tools in the estimations.
For the nonnegative functions $\hat{f},\hat{g},\hat{h}$, we define
\begin{equation}
I^s_k(\hat{f},\hat{g},\hat{h})=\displaystyle\int_D
m_k(\xi,\xi_1,\xi_2)^s
\hat{f}(\xi,\tau)\hat{g}(\xi_1,\tau_1)\hat{h}(\xi_2,\tau_2)
\label{Isk}
\end{equation}
for  $k=1,\cdots, 7$, where the set
$D=\{(\xi_1,\xi_2,\tau_1,\tau_2):\xi=\xi_1+\xi_2,
\tau=\tau_1+\tau_2\}$, and the multipliers $m_k$ are defined as
$$
\begin{array}{c}
m_1=|3\xi_2^2+2\xi_1|; \quad m_2=|3\xi_2^2+2\xi|; \quad m_3=|\xi_2|;\\
m_4=|\xi|; \quad m_5=|3\xi^2-2\xi_2|; \quad m_6=|3\xi^2+2\xi_1|;\\
m_7=|\xi_1-\xi_2||\xi_1+\xi_2|.
\end{array}
$$
We state the relevant estimates on them.
\begin{lem}
Let $f,g,h$ are reasonable functions, then
$$
\begin{array}{c}
(1)\quad I^{\frac{1}{2}}_1(\hat{f},\hat{g},\hat{h})
   \lesssim
              \|f\|_{L^2}\>
              \|g\|_{{X}_{0,\frac{1}{2}+}}\>
              \|h\|_{{Y}_{0,\frac{1}{2}+}};\quad
I^{\frac{1}{2}}_2(\hat{f},\hat{g},\hat{h})
  \lesssim
              \|f\|_{{X}_{0,\frac{1}{2}+}}\>
              \|g\|_{L^2}\>
              \|h\|_{{Y}_{0,\frac{1}{2}+}}; \\
I^{\frac{1}{2}}_3(\hat{f},\hat{g},\hat{h})
   \lesssim
              \|f\|_{{X}_{0,\frac{1}{2}+}}\>
              \|g\|_{{X}_{0,\frac{1}{2}+}}\>
              \|h\|_{L^2};\\
(2)\quad I^{\frac{1}{2}}_4(\hat{f},\hat{g},\hat{h})
   \lesssim
              \|f\|_{L^2}\>
              \|g\|_{{X}_{0,\frac{1}{2}+}}\>
              \|h\|_{{X}_{0,\frac{1}{2}+}^{-}};\quad
I^{\frac{1}{2}}_5(\hat{f},\hat{g},\hat{h})
  \lesssim
              \|f\|_{{Y}_{0,\frac{1}{2}+}}\>
              \|g\|_{L^2}\>
              \|h\|_{{X}_{0,\frac{1}{2}+}^{-}}; \\
I^{\frac{1}{2}}_6(\hat{f},\hat{g},\hat{h})
   \lesssim
              \|f\|_{{Y}_{0,\frac{1}{2}+}}\>
              \|g\|_{{X}_{0,\frac{1}{2}+}}\>
              \|h\|_{L^2};\\
(3)\quad I^{\frac{1}{2}}_7(\hat{f},\hat{g},\hat{h})
   \lesssim
              \|f\|_{L^2}\>
              \|g\|_{{Y}_{0,\frac{1}{2}+}}\>
              \|h\|_{{Y}_{0,\frac{1}{2}+}}.
              \qquad\qquad\qquad\qquad\qquad\qquad\qquad\qquad\qquad\\
\end{array}
$$
\end{lem}
\begin{proof}
We use the argument in \cite{CKS2} to prove the lemma.
For $I^{\frac{1}{2}}_1$, we change variables by setting
$$
\tau=\lambda-\xi^2,\quad \tau_1=\lambda_1-\xi_1^2,\quad
\tau_2=\lambda_2+\xi_2^3,
$$
then,  $I^{\frac{1}{2}}_1(f,g,h)$ is changed into
\begin{equation}
\displaystyle\int m_1^{\frac{1}{2}}\,
\hat{f}(\xi_1+\xi_2,\lambda_1+\lambda_2-\xi_1^2+\xi_2^3)\,
\hat{g}(\xi_1,\lambda_1-\xi_1^2)\,
\hat{h}(\xi_2,\lambda_2+\xi_2^3)\, d \xi_1 d \xi_2 d \lambda_1 d
\lambda_2.\label{2.9}
\end{equation}
We change variables again as follows. Let
\begin{equation}
 (\eta,\omega)=T(\xi_1,\xi_2),\label{transform}
\end{equation}
where
\beqs
\begin{split}
 \eta & = T_1(\xi_1,\xi_2)=\xi_1+\xi_2,\\
 \omega & = T_2(\xi_1,\xi_2)=\lambda_1+\lambda_2-\xi_1^2+\xi_2^3.
\end{split}
\eeqs
Then the Jacobian $J$ of this transform satisfies
$$
|J|=|3\xi_2^2+2\xi_1|.
$$
Define
$$
H(\eta,\omega,\lambda_1,\lambda_2)=
\hat{g}\hat{h}\circ
T^{-1}(\eta,\omega,\lambda_1,\lambda_2),
$$
then, by using $|J|^{\frac{1}{2}}$ to eliminate
$m_1^{\frac{1}{2}}$, (\ref{2.9}) has a bound of
\begin{equation}
    \displaystyle\int\!\! \hat{f}(\eta,\omega)
    \cdot \dfrac{ H(\eta,\omega,\lambda_1,\lambda_2)}
    {|J|^{\frac{1}{2}}}
    \,d\eta d\omega d\lambda_1 d\lambda_2.\label{2.11}
\end{equation}
By H\"{o}lder' inequality, we have
\renewcommand{\arraystretch}{2}
\beqs
\begin{split}
(\ref{2.11})\leq
 &\  \left\|\hat{f}\right\|_{L^2_{\eta\omega}}
     \cdot\displaystyle\int\!\!
      \Big(\int \dfrac{|H(\eta,\omega,\lambda_1,\lambda_2)|^2}{|J|}
      \,d\eta \omega\Big)^{\frac{1}{2}}
      \,d\lambda_1 d\lambda_2\\
=
 &\   \left\|\hat{f}\right\|_{L^2_{\eta\omega}}
      \cdot
      \displaystyle\int\!
      \left\|\hat{g}(\xi_1,\lambda_1-\xi_1^2)\right\|_{L^2_{\xi_1}}d\lambda_1
      \cdot
      \displaystyle\int\!
      \left\|\hat{h}(\xi_2,\lambda_2+\xi_2^3)\right\|_{L^2_{\xi_2}}d\lambda_2\\
\lesssim
 &\   \|f\|_{L^2}\>
      \|g\|_{{X}_{0,\frac{1}{2}+}}\>
      \|h\|_{{Y}_{0,\frac{1}{2}+}},
\end{split}
\eeqs
where we employed the inverse transform of (\ref{transform}) in the
second step and H\"{o}lder' inequality in the third step.

For $I^{\frac{1}{2}}_2$, the modification of the proof is
replacing the variable transform $(\eta,\omega)$ by
\beqs
\begin{split}
 \eta & = T_1(\xi,\xi_2)=\xi-\xi_2,\\
 \omega & = T_2(\xi,\xi_2)=\lambda-\lambda_2-\xi^2-\xi_2^3.
\end{split}
\eeqs
Then the Jacobian $J$ in this situation satisfies
$$
|J|=|3\xi_2^2+2\xi|.
$$
Therefore, we have the claim by the same argument as above but separating
$\hat{g}$ from the integration in this time.

For $I^{\frac{1}{2}}_3$, we take $(\eta,\omega)$ in this time that
\beqs
\begin{split}
 \eta & = T_1(\xi,\xi_1)=\xi-\xi_1,\\
 \omega & = T_2(\xi,\xi_1)=\lambda-\lambda_1-\xi^2+\xi_1^2.
\end{split}
\eeqs
Then the Jacobian $J$ in this situation satisfies
$$
|J|=2|\xi_2|.
$$
So the claim follows again.

For $I^{\frac{1}{2}}_4,I^{\frac{1}{2}}_5,I^{\frac{1}{2}}_6$,
we change variables first by setting
\begin{equation}
\tau=\lambda+\xi^3,\quad \tau_1=\lambda_1-\xi_1^2,\quad
\tau_2=\lambda_2+\xi_2^2,\label{lambda}
\end{equation}
then for $I^{\frac{1}{2}}_4$, we change variables again as
\beqs
\begin{split}
 \eta & = T_1(\xi_1,\xi_2)=\xi_1+\xi_2,\\
 \omega & = T_2(\xi_1,\xi_2)=\lambda_1+\lambda_2-\xi_1^2+\xi_2^2.
\end{split}
\eeqs
Then the Jacobian $J$ satisfies
$$
|J|=2|\xi|.
$$

For $I^{\frac{1}{2}}_5$, we change variables as
\beqs
\begin{split}
 \eta & = T_1(\xi,\xi_2)=\xi-\xi_2,\\
 \omega & = T_2(\xi,\xi_2)=\lambda-\lambda_2+\xi^3-\xi_2^2.
\end{split}
\eeqs
Then the Jacobian $J$ satisfies
$$
|J|=|3\xi^2-2\xi_2|.
$$

For $I^{\frac{1}{2}}_6$, we change variables as
\beqs
\begin{split}
 \eta & = T_1(\xi,\xi_1)=\xi-\xi_1,\\
 \omega & = T_2(\xi,\xi_1)=\lambda-\lambda_1+\xi^3+\xi_1^2.
\end{split}
\eeqs
Then the Jacobian $J$ satisfies
$$
|J|=|3\xi^2+2\xi_1|.
$$
Therefore, we have the conclusions in the second term.

For $I^{\frac{1}{2}}_7$, we change variables first by setting
\begin{equation}
\tau=\lambda+\xi^3,\quad \tau_1=\lambda_1+\xi_1^3,\quad
\tau_2=\lambda_2+\xi_2^3,\label{lambda}
\end{equation}
then we change variables again as
\beqs
\begin{split}
 \eta & = T_1(\xi_1,\xi_2)=\xi_1+\xi_2,\\
 \omega & = T_2(\xi_1,\xi_2)=\lambda_1+\lambda_2+\xi_1^3+\xi_2^3.
\end{split}
\eeqs
Then the Jacobian $J$ satisfies
$$
|J|=3|\xi_1^2-\xi_2^2|.
$$
Thus the claim follows by the same argument.
\hfill$\Box$
\end{proof}

When $s=0$, by Lemma 2.1 we have
\begin{eqnarray*}
I^{0}_2(\hat{f},\hat{g},\hat{h})
  &\lesssim
    \|f\|_{L^{8/3}_{xt}}\>\|g\|_{L^2_{xt}}\>\|h\|_{L^{8}_{xt}}
  &\lesssim
    \|f\|_{X_{0,\frac{3}{16}+}}\>\|g\|_{L^2_{xt}}\>\|h\|_{Y_{0,\frac{1}{2}+}};\\
I^{0}_3(\hat{f},\hat{g},\hat{h})
  &\lesssim
    \|f\|_{L^{3}_{xt}}\>\|g\|_{L^6_{xt}}\>\|h\|_{L^{2}_{xt}}
  &\lesssim
    \|f\|_{X_{0,\frac{1}{4}+}}\>\|g\|_{X_{0,\frac{1}{2}+}}\>\|h\|_{L^{2}_{xt}};\\
I^{0}_5(\hat{f},\hat{g},\hat{h})
  &\lesssim
    \|f\|_{L^{3}_{xt}}\>\|g\|_{L^2_{xt}}\>\|h\|_{L^{6}_{xt}}
  &\lesssim
    \|f\|_{Y_{0,\frac{2}{9}+}}\>\|g\|_{L^2_{xt}}\>\|h\|_{X_{0,\frac{1}{2}+}^-};\\
I^{0}_6(\hat{f},\hat{g},\hat{h})
  &\lesssim
    \|f\|_{L^{3}_{xt}}\>\|g\|_{L^6_{xt}}\>\|h\|_{L^{2}_{xt}}
  &\lesssim
    \|f\|_{Y_{0,\frac{2}{9}+}}\>\|g\|_{X_{0,\frac{1}{2}+}}\>\|h\|_{L^{2}_{xt}}.
\end{eqnarray*}
Interpolation between them and the results in Lemma 2.2, we have the following lemma.
\begin{lem}
For any $s\in \left[0,\dfrac{1}{2}\right]$, the following estimates
hold:
$$
\begin{array}{c}
I^{s}_2(\hat{f},\hat{g},\hat{h})
  \lesssim
              \|f\|_{{X}_{0,\rho_2+}}\>
              \|g\|_{L^2}\>
              \|h\|_{{Y}_{0,\frac{1}{2}+}}; \\
I^{s}_3(\hat{f},\hat{g},\hat{h})
   \lesssim
              \|f\|_{{X}_{0,\rho_3+}}\>
              \|g\|_{{X}_{0,\frac{1}{2}+}}\>
              \|h\|_{L^2};\\
I^{s}_5(\hat{f},\hat{g},\hat{h})
  \lesssim
              \|f\|_{{Y}_{0,\rho_5+}}\>
              \|g\|_{L^2}\>
              \|h\|_{{X}_{0,\frac{1}{2}+}^{-}}; \\
I^{s}_6(\hat{f},\hat{g},\hat{h})
   \lesssim
              \|f\|_{{Y}_{0,\rho_6+}}\>
              \|g\|_{{X}_{0,\frac{1}{2}+}}\>
              \|h\|_{L^2},\\
\end{array}
$$
where
$\rho_2\geq \dfrac{3}{16}+\dfrac{5}{8}s$, $\rho_3\geq \dfrac{1}{4}+\dfrac{1}{2}s$
and $\rho_5, \rho_6\geq \dfrac{2}{9}+\dfrac{5}{9}s$.
\end{lem}

Further, we denote $\psi(t)$ to be an even smooth characteristic
function of the interval $[-1,1]$, then we have the following
estimates.
\begin{lem} Let $\delta\in (0,1)$, $s\in \R$, then the following
estimates hold: \bi
\item[{\rm(i)}]
           $\|f\|_{C_t^0 (\R;H_x^s)}
            \lesssim\|f\|_{X_{s,b}(\phi)}$,
           $\forall$ $b\in (\dfrac{1}{2},1];$

\item[{\rm(ii)}]
           $\|\psi(t)U_\phi(t)u_0\|_{X_{s,b}(\phi)}
            \lesssim\|u_0\|_{H^s}$,
           $\forall$ $b\in (\dfrac{1}{2},1];$

\item[{\rm(iii)}]
          $\left\|\psi(t)\displaystyle\int^t_0 U_\phi(t-s)F(s)\,ds
          \right\|_{X_{s,b}(\phi)}
          \lesssim
          \left\|F\right\|_{X_{s,b-1}(\phi)}$,
          $\forall$ $b\in (\dfrac{1}{2},1]$;

\item[{\rm(iv)}]
          $\|\psi(t/\delta)f\|_{X_{s,b}(\phi)}
          \lesssim
          \delta^{b'-b}\|f\|_{X_{s,b'}(\phi)}$,
          $\forall$ $0\leq b\leq b'<\dfrac{1}{2}$.
\ei
\end{lem}
\begin{proof}
See cf. \cite{GTV}, \cite{P}, \cite{CKS2} for the proofs.
\hfill$\Box$
\end{proof}

 \vspace{0.3cm}
\section{Multilinear Estimates and Local Well-posedness}

\subsection{Bi- and Trilinear Estimates}

We begin with two well-known estimates.
\begin{lem}(\cite{BOP2})
Let $s\geq 0$, $b=\dfrac{1}{2}+$, then for any
$u_1,u_2,u_3\in X_{s,b}$,
$$
    \|u_1 u_2 \overline{u_3}\|_{L^2_{xt}}
    \lesssim
    \|u_1\|_{X_{s,b}}\>\|u_2\|_{X_{s,b}}\>\|u_3\|_{X_{s,b}}.
$$
\end{lem}
\begin{lem}(\cite{KPV1})
Let $l>-\dfrac{3}{4}$, $c,c'=\dfrac{1}{2}+$, then for any
$v_1,v_2\in Y_{l,c}$,
$$
    \|\partial_x(v_1 v_2)\|_{Y_{l,c'-1}}
    \lesssim
    \|v_1\|_{Y_{l,c}}\>\|v_2\|_{Y_{l,c}}.
$$
\end{lem}

Now we turn to prove other two bilinear estimates which improve the
results in \cite{BOP}, \cite{CL}.
Before stating the next lemma, we note an arithmetic fact that
$$
(\tau+\xi^2)-(\tau_1+\xi_1^2)-(\tau_2-\xi_2^3)=\xi^2-\xi_1^2+\xi_2^3
=\xi_2(\xi_2^2+\xi_2+2\xi_1),
$$
if $\xi=\xi_1+\xi_2, \tau=\tau_1+\tau_2$.
It implies that one of the following three cases always occur:
\begin{eqnarray}
&&(a)\,|\tau+\xi^2|\gtrsim |\xi_2||\xi_2^2+\xi_2+2\xi_1|;\qquad
(b)\,|\tau_1+\xi_1^2|\gtrsim |\xi_2||\xi_2^2+\xi_2+2\xi_1|;\nonumber\\
&&(c)\,|\tau_2-\xi_2^3|\gtrsim
|\xi_2||\xi_2^2+\xi_2+2\xi_1|.\label{abc}
\end{eqnarray}
\begin{lem}
Let $l\geq -1$, $s-l\leq 2$ when $s\geq 0$, $s+l\geq -2$ when $s<0$,
and $b,b',c=\dfrac{1}{2}+$, then for any
$u\in X_{s,b}$, $v\in Y_{l,c}$,
$$
    \|uv\|_{X_{s,b'-1}}
    \lesssim
    \|u\|_{X_{s,b}}\>\|v\|_{Y_{l,c}}.
$$
\end{lem}
\begin{proof}
By duality and Plancherel's identity, it suffices to show that for any
$f\in X_{0,1-b'}$,
$$
    \displaystyle\int_{D}\frac{\langle\xi\rangle^s}
    {\langle\xi_1\rangle^s\langle\xi_2\rangle^l}\>
    \hat{f}(\xi,\tau)\>\hat{g}(\xi_1,\tau_1)\>\hat{h}(\xi_2,\tau_2)
    \lesssim
    \|f\|_{X_{0,1-b'}}\>\|g\|_{X_{0,b}}\>\|h\|_{Y_{0,c}}\equiv RHS,
$$
where the set $D=\{(\xi_1,\xi_2,\tau_1,\tau_2):\xi=\xi_1+\xi_2,
\tau=\tau_1+\tau_2\}$.
We may assume that $\hat{f},\hat{g},\hat{h}$ are nonnegative, otherwise
we can replace them by their absolute value without lose of generality.
This point of view will also be used at the following
multilinear estimates without any mentioned.
We divide the integral domain $D$ into
three parts by writing
$$
\displaystyle\int_{D}=\displaystyle\int_{D_1}+
\displaystyle\int_{D_2}+ \displaystyle\int_{D_3},
$$
where
$$
\begin{array}{l}
D_1= \{(\xi_1,\xi_2,\tau_1,\tau_2)\in D:
|\xi|,|\xi_1|, |\xi_2|\lesssim 1\},\\
D_2= \{(\xi_1,\xi_2,\tau_1,\tau_2)\in D:
|\xi|\lesssim|\xi_1|, |\xi_1|\gg 1\},\\
D_3= \{(\xi_1,\xi_2,\tau_1,\tau_2)\in D: |\xi|\gg|\xi_1|, |\xi|\gg
1\}.
\end{array}
$$

\noindent\textbf{Estimate in $D_1$. } By Lemma 2.1 we have
\begin{eqnarray*}
  \displaystyle\int_{D_1}
   \sim
    \displaystyle\int_{D_1}
    \hat{f}(\xi,\tau)\>\hat{g}(\xi_1,\tau_1)\>\hat{h}(\xi_2,\tau_2)
   \lesssim
    \|f\|_{L^2_{xt}}\> \|g\|_{L^3_{xt}}\>\|h\|_{L^6_{xt}}
  \lesssim
    RHS.
\end{eqnarray*}

\noindent\textbf{Estimate in $D_2$. } We shall split it into two
cases,
$$
(1):\ s\geq 0;\qquad (2):\ s<0.
$$

For \textbf{(1): $s \geq 0$}, then
\begin{equation}\label{331}
\displaystyle\int_{D_2} \lesssim \displaystyle\int_{D_2}
\langle\xi_2\rangle^{-l}
\hat{f}(\xi,\tau)\>\hat{g}(\xi_1,\tau_1)\>\hat{h}(\xi_2,\tau_2).
\end{equation}
We may assume that $|\xi_2|\geq 1$, otherwise it can be gotten as
$\displaystyle\int_{D_1}$. Thus,
$$
  (\ref{331})
    \sim  \displaystyle\int_{D_2}|\xi_2|^{-l}
      \hat{f}(\xi,\tau)\>\hat{g}(\xi_1,\tau_1)\>\hat{h}(\xi_2,\tau_2).
$$
We divide $D_2$ again into two subregions:
$$
\begin{array}{l}
D_{21}=\{(\xi_1,\xi_2,\tau_1,\tau_2)\in D_2:
|\xi_2^2+\xi_2+2\xi_1|\ll|\xi_2|^2\},\\
D_{22}=\{(\xi_1,\xi_2,\tau_1,\tau_2)\in D_2:
|\xi_2^2+\xi_2+2\xi_1|\gtrsim|\xi_2|^2\}.
\end{array}
$$

\noindent\textbf{Estimate in $D_{21}$. } Note that
$|3\xi_2^2+2\xi_1|\sim |\xi_2|^2$ in $D_{21}$, therefore,
$\displaystyle\int_{D_{21}}$ is equivalent to
\begin{eqnarray*}
  && \displaystyle\int_{D_{21}}|\xi_2|^{-l-1}\cdot
      m_1^{\frac{1}{2}}\>
      \hat{f}(\xi,\tau)\>\hat{g}(\xi_1,\tau_1)\>\hat{h}(\xi_2,\tau_2)\\
  & \lesssim &
      I_1^{\frac{1}{2}}(\hat{f},\hat{g},\hat{h})
    \lesssim
      \|f\|_{L^2}\>\|g\|_{X_{0,\frac{1}{2}+}}\>\|h\|_{Y_{0,\frac{1}{2}+}}
    \lesssim RHS,
\end{eqnarray*}
where we note that $l\geq -1$ in the second step.

\noindent\textbf{Estimate in $D_{22}$. }  By (\ref{abc}), we can
divide the integral domain into three parts again. But they are
similar to each other, we just take (a): $|\tau+\xi^2|\gtrsim
|\xi_2||\xi_2^2+\xi_2+2\xi_1|$ for example, then since $l\geq -1$
and by Lemma 2.1, $\displaystyle\int_{D_{22}}$ is equivalent to
\begin{eqnarray*}
   \displaystyle\int_{D_{22}}
      \langle\tau+\xi^2\rangle^{\frac{1}{3}}
      \hat{f}(\xi,\tau)\>\hat{g}(\xi_1,\tau_1)\>\hat{h}(\xi_2,\tau_2)
   \lesssim
      \|f\|_{X_{0,\frac{1}{3}}}\>\|g\|_{L^4_{xt}}\>\|h\|_{L^4_{xt}}
    \lesssim RHS.
\end{eqnarray*}

Now we turn to consider the case \textbf{(2): $s<0$}. We also
further split $D_2$ into following two parts $D_{21}^\prime$ and
$D_{22}^\prime$,
$$
\begin{array}{l}
D_{21}^\prime=\{(\xi_1,\xi_2,\tau_1,\tau_2)\in D_2:
|\xi|\sim |\xi_1|\},\\
D_{22}^\prime=\{(\xi_1,\xi_2,\tau_1,\tau_2)\in D_2: |\xi|\ll
|\xi_1|\sim|\xi_2|\}.
\end{array}
$$
\noindent\textbf{Estimate in $D_{21}'$. }  We can give the claim as
$s\geq 0$.

\noindent\textbf{Estimate in $D_{22}'$. } We note that
$|\xi_2^2+\xi_2+2\xi_1|\sim |\xi_2|^2$, $m_1\sim m_2\sim |\xi_2|^2$.
By (\ref{abc}), we split $D_{22}'$ into three parts again, but
similarly we only consider (c): $|\tau_2-\xi_2^3|\gtrsim
|\xi_2||\xi_2^2+\xi_2+2\xi_1|$, then by Lemma 2.3,
\begin{eqnarray*}
  \displaystyle\int_{D_{22}^\prime}
  &\sim & \displaystyle\int_{D_{22}^\prime}
      \langle\xi\rangle^s |\xi_2|^{-s-l}
      \hat{f}(\xi,\tau)\>\hat{g}(\xi_1,\tau_1)\>\hat{h}(\xi_2,\tau_2)\\
  & \lesssim & \displaystyle\int_{D_{22}^\prime}
      |\xi_2|^{-s-l}
      \hat{f}(\xi,\tau)\>\hat{g}(\xi_1,\tau_1)\>\hat{h}(\xi_2,\tau_2)\\
  & \lesssim & \displaystyle\int_{D_{22}^\prime}
      |\xi_2|^{-s-l-3c-s_1}\> m_3^{s_1}
      \hat{f}(\xi,\tau)\>\hat{g}(\xi_1,\tau_1)\>
      \langle\tau_2-\xi_2^3\rangle^c \hat{h}(\xi_2,\tau_2)\\
  & \lesssim &
      I_3^{s_1}\left(\hat{f},\hat{g},\langle\tau-\xi^3\rangle^c
      \hat{h}\right)\\
  & \lesssim & RHS,
\end{eqnarray*}
where $s_1=\dfrac{1}{2}-$, such that $1-b'>
\dfrac{1}{4}+\dfrac{1}{2}s_1$ and $s+l\geq -3c-s_1$ (ensured by $s+l\geq-2$).

\noindent\textbf{Estimate in $D_{3}$. } We have,
$$
\displaystyle\int_{D_3}
  \sim \displaystyle\int_{D_{3}}
      \langle\xi_1\rangle^{-s} |\xi_2|^{s-l}
      \hat{f}(\xi,\tau)\>\hat{g}(\xi_1,\tau_1)\>\hat{h}(\xi_2,\tau_2).
$$
It is much similar to $\displaystyle\int_{D_{22}^\prime}$. Therefore,
when $s\geq 0$, then
$$
\displaystyle\int_{D_{3}}\lesssim  |\xi_2|^{s-l}
\hat{f}(\xi,\tau)\>\hat{g}(\xi_1,\tau_1)\>\hat{h}(\xi_2,\tau_2),
$$
thus we have the claim by
noting that $s-l\leq 2$;
when $s< 0$, since $|\xi_1|\lesssim |\xi_2|$, we have
$$
\displaystyle\int_{D_{3}}
  \lesssim  \displaystyle\int_{D_{3}}
      |\xi_2|^{-l}
      \hat{f}(\xi,\tau)\>\hat{g}(\xi_1,\tau_1)\>\hat{h}(\xi_2,\tau_2),
$$
thus we have the claim again by noting $l\geq -2$.
\hfill$\Box$
\end{proof}

In the next proof of lemma we will use the following algebraic relation
$$
(\tau-\xi^3)-(\tau_1+\xi_1^2)-(\tau_2-\xi_2^2)=-\xi^3-\xi_1^2+\xi_2^2
=-\xi(\xi^2+\xi-2\xi_2),
$$
if $\xi=\xi_1+\xi_2, \tau=\tau_1+\tau_2$.
It implies that one of the following three cases must occur:
\begin{eqnarray}
&&(a)\,|\tau-\xi^3| \gtrsim  |\xi||\xi^2+\xi-2\xi_2|;\qquad
(b)\,|\tau_1+\xi_1^2|\gtrsim |\xi||\xi^2+\xi-2\xi_2|;\nonumber\\
&&(c)\,|\tau_2-\xi_2^2|\gtrsim |\xi||\xi^2+\xi-2\xi_2|.\label{abc2}
\end{eqnarray}
\begin{lem}
Let $s> -\dfrac{1}{4}$, $l<4s$, $l-s< 1$ when $s\geq 0$, $l-2s<1$
when $s<0$, and $b,c'=\dfrac{1}{2}+$, then for any $u_1,u_2\in X_{s,b}$,
$$
    \|\partial_x(u_1\overline{u_2})\|_{Y_{l,c'-1}}
    \lesssim
    \|u_1\|_{X_{s,b}}\>\|u_2\|_{X_{s,b}}.
$$
\end{lem}
\begin{proof}By duality and Plancherel's identity and the fact
$
\widehat{\bar{u}}(\xi,\tau)=\overline{\hat{u}}(-\xi,-\tau),
$
it suffices to show that
$$
    \displaystyle\int_{D}\frac{|\xi|\langle\xi\rangle^l}
    {\langle\xi_1\rangle^s\langle\xi_2\rangle^s}\>
    \hat{f}(\xi,\tau)\>\hat{g}(\xi_1,\tau_1)\>\hat{h}(\xi_2,\tau_2)
    \lesssim
    \|f\|_{Y_{0,1-c'}}\>\|g\|_{X_{0,b}}\>\|h\|_{X_{0,b}^-}\equiv RHS,
$$
where the set $D=\{(\xi_1,\xi_2,\tau_1,\tau_2):\xi=\xi_1+\xi_2,
\tau=\tau_1+\tau_2\}$. We may assume that $|\xi_1|\geq |\xi_2|$ (the
other is similar). We divide the integration domain $D$ into three
parts:
$$
\begin{array}{l}
D_1= \{(\xi_1,\xi_2,\tau_1,\tau_2)\in D:
|\xi|,|\xi_1|, |\xi_2|\lesssim 1\},\\
D_2= \{(\xi_1,\xi_2,\tau_1,\tau_2)\in D:
|\xi|\ll|\xi_1|\sim |\xi_2|, |\xi_1|\gg 1\},\\
D_3= \{(\xi_1,\xi_2,\tau_1,\tau_2)\in D: |\xi|\sim|\xi_1|,
|\xi_1|\gg 1\}.
\end{array}
$$

\noindent\textbf{Estimate in $D_{1}$. }  By Lemma 2.1 we have
\begin{eqnarray*}
  \displaystyle\int_{D_1}
   \sim
    \displaystyle\int_{D_1}
    \hat{f}(\xi,\tau)\>\hat{g}(\xi_1,\tau_1)\>\hat{h}(\xi_2,\tau_2)
   \lesssim
    \|f\|_{L^4_{xt}}\> \|g\|_{L^4_{xt}}\>\|h\|_{L^2_{xt}}
  \lesssim
    RHS.
\end{eqnarray*}

\noindent\textbf{Estimate in $D_{2}$. }  We have
$$
  \displaystyle\int_{D_2}
  \sim
  \displaystyle\int_{D_2}
  |\xi|\langle\xi\rangle^{l}|\xi_2|^{-2s}
  \hat{f}(\xi,\tau)\>\hat{g}(\xi_1,\tau_1)\>\hat{h}(\xi_2,\tau_2)
$$
We divide $D_2$ again into two subparts:
$$
\begin{array}{l}
D_{21}=\{(\xi_1,\xi_2,\tau_1,\tau_2)\in D_2:
|\xi^2+\xi-2\xi_2|\ll|\xi_2|\},\\
D_{22}=\{(\xi_1,\xi_2,\tau_1,\tau_2)\in D_2:
|\xi^2+\xi-2\xi_2|\gtrsim|\xi_2|\}.
\end{array}
$$

\noindent\textbf{Estimate in $D_{21}$. } We note that $|\xi|^2\sim
|\xi_2|$ and $|3\xi^2-2\xi_2|\sim |\xi_2|$ in $D_{21}$, therefore,
\begin{eqnarray*}
      \int_{D_{21}}
  &\sim &
      \displaystyle\int_{D_{21}}|\xi_2|^{-2s+\frac{l+1}{2}}\>
      \hat{f}(\xi,\tau)\>\hat{g}(\xi_1,\tau_1)\>\hat{h}(\xi_2,\tau_2)\\
  &\lesssim&
      \displaystyle\int_{D_{21}}|\xi_2|^{-2s+\frac{l+1}{2}-s_2}\cdot
      m_5^{s_2}\>
      \hat{f}(\xi,\tau)\>\hat{g}(\xi_1,\tau_1)\>\hat{h}(\xi_2,\tau_2)\\
  & \lesssim &
      I_5^{s_2}(\hat{f},\hat{g},\hat{h})\\
  & \lesssim &
      \|f\|_{Y_{0,1-c'}}\>\|g\|_{L^2}\>\|h\|_{X_{0,\frac{1}{2}+}^-}\\
  &  \lesssim &
  RHS,
\end{eqnarray*}
where we use Lemma 2.3 in the fourth step and choose
$s_2=\dfrac{1}{2}-$, such that $1-c'>
\dfrac{2}{9}+\dfrac{5}{9}s_2$ and $-2s+\dfrac{l+1}{2}-s_2\geq 0$ (which
is ensured by $l<4s$).

\noindent\textbf{Estimate in $D_{22}$. } In $D_{22}$, all the cases
of $|\xi|^2\ll |\xi_2|$, $|\xi|^2\sim |\xi_2|$ or $|\xi|^2\gg
|\xi_2|$, and (a), (b) or (c) in (\ref{abc2}) maybe occur. Note that
$ |\xi^2+\xi-2\xi_2|\sim \max\{|\xi|^2, |\xi_2|\}$, we see that the
worst case is: $|\xi_2|\gg |\xi|^2$ and $|\tau-\xi^3|\gg
\max\{|\tau_1+\xi_1^2|, |\tau_2-\xi_2^2|\}$. We only consider the
integration under this part (the others can be treated similarly but
employing the estimates on $I^s_5$ or $I^s_6$ in Lemma 2.3). Then
$|\tau-\xi^3|\gtrsim |\xi||\xi_2|$, and thus,
\begin{eqnarray*}
\int_{D_{22}}
  &\sim&
      \displaystyle\int_{D_{22}}|\xi|^{c'}\langle\xi\rangle^l
      |\xi_2|^{-2s-(1-c')}\>\langle\tau-\xi^3\rangle^{1-c'}
      \hat{f}(\xi,\tau)\>\hat{g}(\xi_1,\tau_1)\>\hat{h}(\xi_2,\tau_2)\\
  &\lesssim&
      \displaystyle\int_{D_{22}}
      \langle\tau-\xi^3\rangle^{1-c'}\>
      \hat{f}(\xi,\tau)\>\hat{g}(\xi_1,\tau_1)\>\hat{h}(\xi_2,\tau_2)\\
  & \lesssim &
      \|f\|_{Y_{0,1-c'}}\>\|g\|_{L^4_{xt}}\>\|h\|_{L^4_{xt}}\\
  &  \lesssim &
  RHS,
\end{eqnarray*}
where we note $s>-\dfrac{1}{4}, l<4s$ in the second step and use
Lemma 2.1 in the fourth step.

\noindent\textbf{Estimate in $D_{3}$. }  We have
\begin{equation}\label{341}
    \int_{D_{3}}\sim \int_{D_{3}}|\xi|^{1+l-s}\langle\xi_2\rangle^{-s}
    \hat{f}(\xi,\tau)\>\hat{g}(\xi_1,\tau_1)\>\hat{h}(\xi_2,\tau_2).
\end{equation}
We split it into two cases to analysis,
$$
(1):\ s\geq 0; \quad (2):\ s< 0.
$$

For \textbf{(1): $s\geq 0$}, then
\begin{equation}\label{342}
    (\ref{341})\lesssim\displaystyle\int_{D_{3}}|\xi|^{1+l-s}
    \hat{f}(\xi,\tau)\>\hat{g}(\xi_1,\tau_1)\>\hat{h}(\xi_2,\tau_2).
\end{equation}
In $D_3$, we have $|\xi^2+\xi-2\xi_2|\sim |\xi|^2$, $m_5\sim m_6\sim
|\xi|^2$. By (\ref{abc2}), we split the domain $D_3$ into three
parts, but we only take (a) for example. Therefore, by Lemma 2.2,
\begin{eqnarray*}
  (\ref{342})&\lesssim&
      \displaystyle\int_{D_{3}}|\xi|^{1+l-s+3c'-3-\frac{1}{2}}
      \>m_4^{\frac{1}{2}}\>
      \langle\tau-\xi^3\rangle^{1-c'}\hat{f}(\xi,\tau)\>
      \hat{g}(\xi_1,\tau_1)\>\hat{h}(\xi_2,\tau_2)\\
  &\lesssim &
      I_4^{\frac{1}{2}}(\langle\tau-\xi^3\rangle^{1-c'}\hat{f},
      \hat{g},\hat{h}) \\
  & \lesssim &
      RHS,
\end{eqnarray*}
where we note that $l-s<1$ and $c'=\dfrac{1}{2}+$, thus
$l-s+3c'-\dfrac{5}{2}\leq 0$.

For \textbf{(2): $s< 0$}, since $|\xi_2|\lesssim |\xi|$, we have,
$$
    (\ref{341})\lesssim\displaystyle\int_{D_{3}}|\xi|^{1+l-2s}
    \hat{f}(\xi,\tau)\>\hat{g}(\xi_1,\tau_1)\>\hat{h}(\xi_2,\tau_2).
$$
Similar to (1), we have the claim since $l-2s<1$.
 \hfill$\Box$
\end{proof}

\subsection{Local Well-posedness}
Write the sets
\begin{eqnarray*}
  R_{\beta} &\equiv& \{(s,l):s\geq 0\};\quad
  R_{kdv} \equiv  \left\{(s,l):l> -\frac{3}{4}\right\};  \\
  R_{\alpha} &\equiv& \{(s,l):l\geq -1; s-l\leq 2 \mbox{\,\,when\,\,} s\geq 0,
                   s+l\geq -2 \mbox{\,\,when\,\,} s< 0\};   \\
  R_{\gamma} &\equiv& \left\{(s,l):s> -\frac{1}{4}; l<4s; l-s< 1
                   \mbox{\,\,when\,\,} s\geq 0,
                   l-2s<1 \mbox{\,\,when\,\,} s< 0\right\},
\end{eqnarray*}
and let
\begin{eqnarray*}
  R_{\beta=0}&\equiv& R_{kdv}\cap R_{\alpha}\cap R_{\gamma}
    =\left\{(s,l):l>-\frac{3}{4}, l<4s, s-2\leq l < s+1\right\};\\
  R_{\beta\neq0}&\equiv& R_{kdv}\cap R_{\alpha}\cap R_{\gamma}\cap
  R_{\beta}
    =\left\{(s,l):l>-\frac{3}{4}, s\geq 0, l<4s, s-2\leq l < s+1\right\}.
\end{eqnarray*}
We assume that $(s,l)\in R_{\beta=0}$ in the case of $\beta=0$ and
$(s,l)\in R_{\beta\neq0}$ in the case of $\beta\neq0$.  Define the
maps
\begin{eqnarray*}
  \Phi_1(u,v) &=& \psi(t)S(t)u_0-i\psi(t)\displaystyle\int_0^t
  S(t-t')\psi(t'/\delta)\left[\alpha (uv)(t')+\beta(|u|^2u)(t')\right]\,dt', \\
  \Phi_2(u,v) &=& \psi(t)W(t)v_0+\psi(t)\displaystyle\int_0^t
  W(t-t')\psi(t'/\delta)\left[\gamma \partial_x(|u|^2)(t')
  -\frac{1}{2}\partial_x(v^2)(t')\right]\,dt',
\end{eqnarray*}
where $S(t), W(t)$ are $U_\phi(t)$ with $\phi=\xi^2,-\xi^3$
respectively, then (\ref{NLS-KdV}) is locally well-posed on
$[-\delta, \delta]$ only if $(\Phi_1,\Phi_2)$ has a unique fixed
point. By Lemma 2.4, Lemmas 3.1--3,4, the fixed point theory and a
standard process (see \cite{CL} cf.), we prove Theorem 1.1 with the
estimates on the lifetime and solutions that for
$\mu=\max\{b'-b,c'-c\}>0$,
$$
\delta\sim (\|u_0\|_{H^s}+\|v_0\|_{H^l})^{-\mu};\quad
\|u\|_{X_{s,\frac{1}{2}+}^\delta}+\|v\|_{Y_{l,\frac{1}{2}+}^\delta}\lesssim
\|u_0\|_{H^s}+\|v_0\|_{H^l},
$$
or
$$
\delta\sim
\min\left\{\|u_0\|_{H^s},\|v_0\|_{H^l},
\frac{\|v_0\|_{H^l}^2}{\|u_0\|_{H^s}}\right\}^{-\mu};\
\|u\|_{X_{s,\frac{1}{2}+}^\delta}\lesssim \|u_0\|_{H^s},
\|v\|_{Y_{l,\frac{1}{2}+}^\delta}\lesssim \|v_0\|_{H^l}.
$$

\noindent{\it Remark.} From the proof of Lemma 3.3 and Lemma 3.4,
more general conditions for the local well-posedness
region (see Figure 1) on the top-left and bottom-right areas are
$$
3c+\frac{3}{2}-2b'<l-s\leq \frac{5}{2}-3c',
$$
for any $c,c',b'$ large and suitable close to $\dfrac{1}{2}$, and $c'>c$.
However, one always has the restriction that
$$
\mbox{ if } b< l-s\leq a,\quad  \mbox{ then } a+b\leq 3.
$$
This implies that the well-posedness region is contained in a belt with
the distance of 3.
 \vspace{0.3cm}
\section{The Proof of Theorem 1.3}

In this section, we consider the global well-posedness of the
solutions obtained in Theorem 1.1 when $\alpha\gamma>0$ and
$(u_0,v_0)\in H^s(\R)\times H^s(\R)$ for some $s>0$. We assume that
$v$ is real valued from now.
In this paper, we pursue to lengthen the
lifetime by some special techniques and a useful
conservation law on $L^2$-norm of $u$. We do nothing on the almost
conserved quantities but cite what obtained in \cite{P} directly.

\subsection{Some Variant Multilinear Estimates}

Now we turn to establish some special multilinear estimates, which
are useful in the next subsection although there is a bit cumbersome in
some.

We will use the following two inequalities frequently in the
multilinear estimates below, which follow from Lemma 2.4 (iv) and
Lemma 6.3. They are, \beq
 \|f\|_{X_{s,b}^\delta(\phi)}\lesssim \delta^{(\frac{1}{2}-b)-}
 \|f\|_{X_{s,\frac{1}{2}+}^\delta(\phi)};\quad
 \|\psi(t/\delta)f\|_{X_{s,b}(\phi)}
 \lesssim
 \delta^{b'-b}\|f\|_{X_{s,b'}(\phi)}\label{L4.2}
\eeq
for $b,b'\in [0,\dfrac{1}{2})$ with $b'\geq b$.
\begin{lem}
Let $s\geq 0$, $c=\dfrac{1}{2}+$, $\delta\in (0,1)$, then for any
$v_1,v_2\in Y_{s,c}^\delta$,
$$
    \|\psi(t/\delta)\partial_x(\widetilde{v_1} \widetilde{v_2})\|_{Y_{s,c-1}}
    \lesssim
    \delta^{\frac{1}{2}-}\cdot\|v_1\|_{Y_{s,c}^\delta}
    \>\|v_2\|_{Y_{s,c}^\delta},
$$
where $\widetilde{v_1}, \widetilde{v_2}$ are the extensions of
$v_1|_{t\in [-\delta,\delta]}$ and $v_2|_{t\in [-\delta,\delta]}$
such that $\|v_1\|_{Y_{s,c}^\delta}=\|\widetilde{v_1}\|_{Y_{s,c}}$,
$\|v_2\|_{Y_{s,c}^\delta}=\|\widetilde{v_2}\|_{Y_{s,c}}$.
\end{lem}
\begin{proof}
By duality and Plancherel's identity,
it suffices to show that for any $f\in Y_{0,1-c}$,
$$
    \displaystyle\int_{D}\frac{|\xi|\langle\xi\rangle^s}
    {\langle\xi_1\rangle^s\langle\xi_2\rangle^s}\>
    \widehat{f_\delta}(\xi,\tau)\>\hat{g}(\xi_1,\tau_1)\>\hat{h}(\xi_2,\tau_2)
    \lesssim
    \delta^{\frac{1}{2}-}\cdot
    \|f\|_{Y_{0,1-c}}\>\|g\|_{Y_{0,c}}\>\|h\|_{Y_{0,c}}\equiv RHS,
$$
where the set $D=\{(\xi_1,\xi_2,\tau_1,\tau_2):\xi=\xi_1+\xi_2,
\tau=\tau_1+\tau_2\}$, $f_\delta=\psi(t/\delta)f$, and
$$
\hat{g}(\xi,\tau)=\langle\xi\rangle^s\widehat{\widetilde{v_1}}(\xi,\tau);\quad
\hat{h}(\xi,\tau)=\langle\xi\rangle^s\widehat{\widetilde{v_2}}(\xi,\tau).
$$
Further, since $s\geq 0$, it is sufficient to show
$$
    \displaystyle\int_{D}|\xi|\>
    \widehat{f_\delta}(\xi,\tau)\>\hat{g}(\xi_1,\tau_1)\>\hat{h}(\xi_2,\tau_2)
    \lesssim
    RHS.
$$
We may assume that $|\xi_1|\geq |\xi_2|$ by
symmetry and write
$$
\displaystyle\int_{D}=\displaystyle\int_{D_1}+ \displaystyle\int_{D_2},
$$
where
$$
D_1= \{(\xi_1,\xi_2,\tau_1,\tau_2)\in D: |\xi|\lesssim|\xi_2|\};\quad
D_2=\{(\xi_1,\xi_2,\tau_1,\tau_2)\in D: |\xi|\ll|\xi_2|\sim|\xi_1|\}.
$$

First, by the arithmetic fact
$$
(\tau-\xi^3)-(\tau_1-\xi_1^3)-(\tau_2-\xi_2^3)
=-3\xi\xi_1\xi_2,
$$
we can split $D_1$ into three parts:
$$
(a)\,|\tau-\xi^3|\gtrsim |\xi||\xi_1||\xi_2|;\ \
(b)\,|\tau_1-\xi_1^3|\gtrsim |\xi||\xi_1||\xi_2|;\ \
(c)\,|\tau_2-\xi_2^3|\gtrsim |\xi||\xi_1||\xi_2|.
$$
We just take (a) for example, since the other two are similar, then
$\displaystyle\int_{D_1}$ is controlled by
\begin{eqnarray*}
  \displaystyle\int_{D_1}\langle\tau-\xi^3\rangle^{\frac{1}{3}}
  \widehat{f_\delta}(\xi,\tau)\>\hat{g}(\xi_1,\tau_1)\>\hat{h}(\xi_2,\tau_2)
 \lesssim
  \|f_\delta\|_{Y_{0,\frac{1}{3}}}\>\|g\|_{L^4_{xt}}\>\|h\|_{L^4_{xt}}
  \qquad\\
 \lesssim
  \|f_\delta\|_{Y_{0,\frac{1}{3}}}\>\|g\|_{Y_{0,\frac{1}{3}+}}
  \>\|h\|_{Y_{0,\frac{1}{3}+}}
 \lesssim
  RHS
\end{eqnarray*}
by Lemma 2.1 and (\ref{L4.2}).
Second, note that $\xi_1\cdot\xi_2<0$, by Lemma 2.2 and (\ref{L4.2}), we have
\begin{eqnarray*}
  \displaystyle\int_{D_2}
 \sim
  I^{\frac{1}{2}}_7(\widehat{f_\delta},\hat{g},\hat{h})
 \lesssim
  \|f_\delta\|_{L^2_{xt}}\>\>\|g\|_{Y_{0,\frac{1}{2}+}}
  \>\|h\|_{Y_{0,\frac{1}{2}+}}
 \lesssim
  RHS.
\end{eqnarray*}
This completes the proof of the lemma.  \hfill$\Box$
\end{proof}

By a general result in \cite{CKSTT3}, Lemma 6.1 and Lemma 4.1, we have
\begin{cor}
Let $I=I_{N,s}$, $s\geq 0$, $c=\dfrac{1}{2}+$, $\delta\in (0,1)$, then
for any $v_1,v_2\in Y_{s,c}^\delta$,
$$
    \|\psi(t/\delta)\partial_xI(\widetilde{v_1} \widetilde{v_2})\|_{Y_{1,c-1}}
    \lesssim
    \delta^{\frac{1}{2}-}\cdot\|Iv_1\|_{Y_{1,c}^\delta}
    \>\|Iv_2\|_{Y_{1,c}^\delta},
$$
where $\widetilde{v_1}, \widetilde{v_2}$ are same with Lemma 4.1.
\end{cor}

\begin{lem}
Let $s\geq 0$, $b,c=\dfrac{1}{2}+$, $\delta\in (0,1)$, then for any
$u\in X_{s,b}^\delta$, $v\in Y_{s,c}^\delta$,
$$
    \|\psi(t/\delta)\>\tilde{u} \tilde{v}\|_{X_{s,b-1}}
    \lesssim
    \delta^{\frac{13}{16}-}\cdot\|u\|_{X_{s,b}^\delta}
    \>\|v\|_{Y_{s,c}^\delta},
$$
where $\tilde{u}, \tilde{v}$ are the extensions of
$u|_{t\in [-\delta,\delta]}$ and $v|_{t\in [-\delta,\delta]}$
such that $\|u\|_{X_{s,b}^\delta}=\|\tilde{u}\|_{X_{s,b}}$,
$\|v\|_{Y_{s,c}^\delta}=\|\tilde{v}\|_{Y_{s,c}}$.
\end{lem}
\begin{proof}
Since $s\geq0$, it suffices to prove that for any $f\in X_{0,1-b}$,
\begin{equation}\label{4.4}
    \displaystyle\int_{D}
    \widehat{f_\delta}(\xi,\tau)\>\hat{g}(\xi_1,\tau_1)\>\hat{h}(\xi_2,\tau_2)
    \lesssim
    \delta^{\frac{13}{16}-}\cdot
    \|f\|_{X_{0,1-b}}\>\|g\|_{X_{0,b}}\>\|h\|_{Y_{0,c}}\equiv RHS,
\end{equation}
where the set $D=\{(\xi_1,\xi_2,\tau_1,\tau_2):\xi=\xi_1+\xi_2,
\tau=\tau_1+\tau_2\}$, $f_\delta=\psi(t/\delta)f$, and
$$
\hat{g}(\xi,\tau)=\langle\xi\rangle^s\widehat{\tilde{u}}(\xi,\tau);\quad
\hat{h}(\xi,\tau)=\langle\xi\rangle^s\widehat{\tilde{v}}(\xi,\tau).
$$
But the left hand side of (\ref{4.4}) is controlled by
$$
\|f_\delta\|_{L^2_{xt}}\>\|g\|_{L^{\frac{8}{3}}_{xt}}\>\|h\|_{L^8_{xt}}
\lesssim
\|f_\delta\|_{X_{0,0}}\>\|g\|_{X_{0,\frac{5}{16}+}}\>\|h\|_{Y_{0,\frac{1}{2}+}}
\lesssim
RHS,
$$
by using Lemma 2.1 and (\ref{L4.2}).  \hfill$\Box$
\end{proof}

Again, by the general result in \cite{CKSTT3}, Lemma 6.1 and Lemma 4.3, we have
\begin{cor}
Let $I=I_{N,s}$, $s\geq 0$, $b,c=\dfrac{1}{2}+$, $\delta\in (0,1)$, then
for any $u\in X_{s,b}^\delta$, $v\in Y_{s,c}^\delta$,
$$
    \left\|\psi(t/\delta)\>I(\tilde{u} \tilde{v})\right\|_{X_{1,b-1}}
    \lesssim
    \delta^{\frac{13}{16}-}\cdot\|Iu\|_{X_{1,b}^\delta}
    \>\|Iv\|_{Y_{1,c}^\delta},
$$
where $\tilde{u}, \tilde{v}$ are same with Lemma 4.3.
\end{cor}

\begin{lem}
Let $I=I_{N,s}$, $s\geq 0$, $b=\dfrac{1}{2}+$, $\delta\in (0,1)$, then
for any $u\in X_{s,b}^\delta$,
$$
    \left\|\psi(t/\delta)\>I(|\tilde{u}|^2 \tilde{u})\right\|_{X_{1,b-1}}
    \lesssim
    \delta^{\frac{1}{2}-}\cdot\|Iu\|_{X_{1,b}^\delta}
    \>\|u\|_{X_{0,b}^\delta}^2,
$$
where $\tilde{u}$ is the extension of
$u|_{t\in [-\delta,\delta]}$
such that $\|u\|_{X_{s,b}^\delta}=\|\tilde{u}\|_{X_{s,b}}$.
\end{lem}
\begin{proof}
By Lemma 2.4 (iv), it suffices to show that
$$
\|I(|u|^2 u)\|_{X_{1,0}}
    \lesssim
    \|Iu\|_{X_{1,b}}
    \>\|u\|_{X_{0,b}}^2
$$
for any $u\in X_{s,b}$. Further, it is equivalent to
\begin{equation}\label{4.5}
\displaystyle\int_D m(\xi)\langle\xi\rangle\>
f(\xi,\tau)\> \hat{u}(\xi_1,\tau_1)\> \hat{u}(\xi_2,\tau_2)
\> \overline{\hat{u}}(-\xi_3,-\tau_3)
\>\lesssim\>
\|f\|_{L^2_{xt}}\>\|Iu\|_{X_{1,b}}\>\|u\|_{X_{0,b}}^2
\end{equation}
for any $f\in L^2(\R^2)$, where
$D=\{(\xi_1,\xi_2,\xi_3,\tau_1,\tau_2,\tau_3):\xi=\xi_1+\xi_2+\xi_3,
\tau=\tau_1+\tau_2+\tau_3\}$. Note that there is at least one of
$\xi_j, j=1,2,3$ such that $|\xi|\lesssim |\xi_j|$.
Without loss of generality, we assume that $j=1$, then
$m(\xi)|\xi|\lesssim m(\xi_1)|\xi_1|$.
Therefore, by Lemma 2.1, (\ref{4.5}) is bounded by
$$
\displaystyle\int_D
f(\xi,\tau)\> \langle\xi_1\rangle\widehat{Iu}(\xi_1,\tau_1)\> \hat{u}(\xi_2,\tau_2)
\>\overline{\hat{u}}(-\xi_3,-\tau_3)
\>\lesssim\>
\|f\|_{L^2_{xt}}\>\|(1+D_x) Iu\|_{L^6_{xt}}\>\|u\|_{L^6_{xt}}^2
\>\lesssim\>
RHS.
$$
This completes the proof of the lemma. \hfill$\Box$
\end{proof}
\begin{lem}
Let $I=I_{N,s}$, $s> 0$, $b,c=\dfrac{1}{2}+$, $\delta\in (0,1)$, then
for any $u\in X_{s,b}^\delta$,
$$
    \left\|\psi(t/\delta)\>\partial_xI(|\tilde{u}|^2 )\right\|_{Y_{1,c-1}}
    \lesssim
    \delta^{a}\cdot\|Iu\|_{X_{1,b}^\delta}
    \>\|u\|_{X_{0,b}^\delta}
    +(N^{-1}\delta^{a}+\delta^{\frac{13}{16}-})
    \cdot\|Iu\|_{X_{1,b}^\delta}^2.
$$
where $a=\left(\dfrac{13}{16}-\dfrac{1}{3}\right)-$, and
$\tilde{u}$ is the extension of
$u|_{t\in [-\delta,\delta]}$
such that $\|u\|_{X_{s,b}^\delta}=\|\tilde{u}\|_{X_{s,b}}$.
\end{lem}
\begin{proof}
By duality and Plancherel's identity, it suffices to show that
for any $f\in Y_{0,1-c}$,
\begin{eqnarray*}
    \displaystyle\int_{D}&&\!\!\!\!\!\!\!\!\!\!|\xi|\langle\xi\rangle m(\xi)\>
    \widehat{f_\delta}(\xi,\tau)\>\widehat{\tilde{u}}(\xi_1,\tau_1)
    \>\overline{\widehat{\tilde{u}}}(-\xi_2,-\tau_2)\\
    &&\lesssim
    \|f\|_{Y_{0,1-c}}
    \left(\delta^{a}\|Iu\|_{X_{1,b}^\delta}
    \|u\|_{X_{0,b}^\delta}+(N^{-1}\delta^{a}
    +\delta^{\frac{13}{16}-})
    \|Iu\|_{X_{1,b}^\delta}^2\right)\equiv RHS,
\end{eqnarray*}
where the set $D=\{(\xi_1,\xi_2,\tau_1,\tau_2):\xi=\xi_1+\xi_2,
\tau=\tau_1+\tau_2\}$, $f_\delta=\psi(t/\delta)f$. We restrict in
$D$ that $|\xi_1|\geq |\xi_2|$ (the other is similar) and divide it
into the following four parts:
\begin{eqnarray*}
D_1&=&\{(\xi_1,\xi_2,\tau_1,\tau_2)\in D:
      |\xi|,|\xi_1|,|\xi_2|\lesssim N\};\\
D_2&=&\{(\xi_1,\xi_2,\tau_1,\tau_2)\in D:
      |\xi|\lesssim N, |\xi_1|\sim|\xi_2|\gg N\};\\
D_3&=&\{(\xi_1,\xi_2,\tau_1,\tau_2)\in D:
      |\xi_2|\lesssim N, |\xi|\sim|\xi_1|\gg N\};\\
D_4&=&\{(\xi_1,\xi_2,\tau_1,\tau_2)\in D:
      |\xi|,|\xi_1|,|\xi_2|\gg N\}.
\end{eqnarray*}

\noindent\textbf{Estimate in $D_{1}$. } Note that
$m(\xi),m(\xi_1),m(\xi_2)\sim 1$ in $D_1$, thus,
$$
\displaystyle\int_{D_1}
\sim
\displaystyle\int_{D_1}|\xi|\langle\xi\rangle\>
\widehat{f_\delta}(\xi,\tau)\>\widehat{\tilde{u}}(\xi_1,\tau_1)
    \>\overline{\widehat{\tilde{u}}}(-\xi_2,-\tau_2).
$$
We divide $D_1$ into two parts again:
\begin{eqnarray*}
D_{11} &=& \{(\xi_1,\xi_2,\tau_1,\tau_2)\in D_1:
      |\xi|\lesssim|\xi_2|\};\\
D_{12} &=& \{(\xi_1,\xi_2,\tau_1,\tau_2)\in D_1:
      |\xi|\gg|\xi_2|\}.
\end{eqnarray*}
\noindent\textbf{Estimate in $D_{11}$. }  By (\ref{XE1}),
(\ref{XE2}) and (\ref{L4.2}), we have
\begin{eqnarray*}
  \displaystyle\int_{D_{11}}
 &\lesssim&
  \displaystyle\int_{D_{11}}\langle\xi_1\rangle\langle\xi_2\rangle\>
  \widehat{f_\delta}(\xi,\tau)\>\widehat{I\tilde{u}}(\xi_1,\tau_1)
    \>\overline{\widehat{I\tilde{u}}}(-\xi_2,-\tau_2)\\
 &\lesssim&
  \|f_\delta\|_{L^p_{xt}}\>
  \|(1+ D_x)I \tilde{u}\|_{L^{p_1}_{xt}}^2
 \lesssim
  \delta^{\frac{13}{16}-}\>\|f\|_{Y_{0,1-c}}\>\|I\tilde{u}\|_{X_{1,b}}^2,
\end{eqnarray*}
where $\dfrac{1}{p}+\dfrac{2}{p_1}=1$ such that
$1-c>\dfrac{4}{3}\left(\dfrac{1}{2}-\dfrac{1}{p}\right)$.

\noindent\textbf{Estimate in $D_{12}$. }   By (\ref{abc2}), we
further split $D_{12}$ into three parts, but each part is similar,
we only take (a) for example, then $|\tau-\xi^3| \gtrsim  |\xi|^3$
and $|\xi|\sim |\xi_1|$, thus we have
\begin{eqnarray*}
  \displaystyle\int_{D_{12}}
 &\sim&
  \displaystyle\int_{D_{12}}\langle\tau-\xi^3\rangle^\frac{1}{3}\>
  \widehat{f_\delta}(\xi,\tau)
    \>\langle\xi_1\rangle\widehat{I\tilde{u}}(\xi_1,\tau_1)
    \>\overline{\widehat{\tilde{u}}}(-\xi_2,-\tau_2)\\
 &\lesssim&
  \left\|\F_{\xi\tau}^{-1}\left(\langle\tau-\xi^3\rangle^\frac{1}{3}
  \widehat{f_\delta}\right)\right\|_{L^q_{xt}}\>
  \|(1+ D_x)I \tilde{u}\|_{L^{q_1}_{xt}}\>\|\tilde{u}\|_{L^{q_1}_{xt}}\\
 &\lesssim&
  \delta^{a}\>\|f\|_{Y_{0,1-c}}
  \>\|I\tilde{u}\|_{X_{1,b}}\>\|\tilde{u}\|_{X_{0,b}},
\end{eqnarray*}
where $\dfrac{1}{q}+\dfrac{2}{q_1}=1$ such that
$1-c-\dfrac{1}{3}>\dfrac{4}{3}\left(\dfrac{1}{2}-\dfrac{1}{q}\right)$.

\noindent\textbf{Estimate in $D_{2}$. }  We will show at the
following that
$$
    \displaystyle\int_{D_2}
    \lesssim
    \delta^{\frac{13}{16}-}\cdot\|f\|_{Y_{0,1-c}}\>\|I\tilde{u}\|_{X_{1,b}}^2.
$$
Indeed, it is sufficient if we show
\begin{equation}\label{4.6}
    \displaystyle\int_{D_2} \frac{\xi\langle\xi\rangle}
    {\langle\xi_1\rangle\langle\xi_2\rangle}\cdot
    \frac{m(\xi)}{m(\xi_1)m(\xi_2)}\cdot
    \widehat{f_\delta}(\xi,\tau)\>\hat{g}(\xi_1,\tau_1)\>\hat{h}(\xi_2,\tau_2)
    \lesssim
    \delta^{\frac{13}{16}-}\>
    \|f\|_{Y_{0,1-c}}\>\|g\|_{X_{0,b}}\>\|h\|_{X_{0,b}^-},
\end{equation}
where $\hat{g}=\widehat{I\tilde{u}}, \hat{h}=\widehat{I\overline{\tilde{u}}}$.
The left hand side of (\ref{4.6}) is equivalent to
\begin{eqnarray*}
 N^{2s-2}\displaystyle\int_{D_2} |\xi|\langle\xi\rangle|\xi_1|^{-2s}
&&\!\!\!\!\!\!\!\!\!\!
 \widehat{f_\delta}(\xi,\tau)\>\hat{g}(\xi_1,\tau_1)\>\hat{h}(\xi_2,\tau_2)
\lesssim
 \displaystyle\int_{D_2}
 \widehat{f_\delta}(\xi,\tau)\>\hat{g}(\xi_1,\tau_1)\>\hat{h}(\xi_2,\tau_2)\\
&&\lesssim
 \delta^{\frac{13}{16}-}\>
 \|f\|_{Y_{0,1-c}}\>\|g\|_{X_{0,b}}\>\|h\|_{X_{0,b}^-},
\end{eqnarray*}
by the same way used in $\displaystyle\int_{D_{11}}$
in the last step.

\noindent\textbf{Estimate in $D_{3}$. }  We have
$$
\int_{D_3}\sim \int_{D_3}|\xi|\>
\widehat{f_\delta}(\xi,\tau)\>\langle\xi_1\rangle\widehat{I\tilde{u}}(\xi_1,\tau_1)
    \>\overline{\widehat{\tilde{u}}}(-\xi_2,-\tau_2).
$$
By the same manner used in $\displaystyle\int_{D_{12}}$, we have
$$
\displaystyle\int_{D_3}
\>\lesssim\>
\delta^{a}\>\|f\|_{Y_{0,1-c}}\>\|I\tilde{u}\|_{X_{1,b}}\>\|\tilde{u}\|_{X_{0,b}}.
$$

\noindent\textbf{Estimate in $D_{4}$. }   We will show in the
following that
$$
    \displaystyle\int_{D_4}\lesssim
    N^{-1}\delta^{a}\cdot\|f\|_{Y_{0,1-c}}\>\|I\tilde{u}\|_{X_{1,b}}^2.
$$
In fact, it suffices to show
$$
    \displaystyle\int_{D_4} \frac{\xi\langle\xi\rangle}
    {\langle\xi_1\rangle\langle\xi_2\rangle}\cdot
    \frac{m(\xi)}{m(\xi_1)m(\xi_2)}\cdot
    \widehat{f_\delta}(\xi,\tau)\>\hat{g}(\xi_1,\tau_1)\>\hat{h}(\xi_2,\tau_2)
    \lesssim
    N^{-1}\delta^{a}\>
    \|f\|_{Y_{0,1-c}}\>\|g\|_{X_{0,b}}\>\|h\|_{X_{0,b}^-}.
$$
The left hand side is controlled by
$$
N^{s-1}\displaystyle\int_{D_4}|\xi|^{s+1}|\xi_1|^{-s}|\xi_2|^{-s}\>
\widehat{f_\delta}(\xi,\tau)\>\hat{g}(\xi_1,\tau_1)\>\hat{h}(\xi_2,\tau_2).
$$
We divide $D_4$ into two subregions again:
$$
\begin{array}{l}
D_{41}=\{(\xi_1,\xi_2,\tau_1,\tau_2)\in D_4:
|\xi^2+\xi-2\xi_2|\ll|\xi|^2\};\\
D_{42}=\{(\xi_1,\xi_2,\tau_1,\tau_2)\in D_4:
|\xi^2+\xi-2\xi_2|\gtrsim|\xi|^2\}.
\end{array}
$$

\noindent\textbf{Estimate in $D_{41}$. }  Note that
$|\xi_2|\sim|\xi|^2$ and $|3\xi^2-\xi_2|\sim|\xi|^2$ in $D_{41}$,
then by Lemma 2.3, we have
\begin{eqnarray*}
  \displaystyle\int_{D_{41}}
 &\lesssim&
  N^{s-1}\displaystyle\int_{D_{41}}|\xi|^{1-3s}
  \widehat{f_\delta}(\xi,\tau)\>\hat{g}(\xi_1,\tau_1)\>\hat{h}(\xi_2,\tau_2)
  \\
 &\lesssim&
  N^{s-1}N^{(-3s)+}I^{\frac{1}{2}-}_5(\widehat{f_\delta},\hat{g},\hat{h}) \\
 &\lesssim&
  N^{(-2s-1)+}\delta^{\frac{1}{2}-}\>
  \|f\|_{Y_{0,1-c}}\>\|g\|_{X_{0,b}}\>\|h\|_{X_{0,b}^-},
\end{eqnarray*}
note that $s>0$ in the second step. Since
$N^{(-2s-1)+}\delta^{\frac{1}{2}-}\leq N^{-1}\delta^{a}$, we obtain
the claim.

\noindent\textbf{Estimate in $D_{42}$. }   By (\ref{abc2}), we can
split $D_{42}$ again into three parts, as above, we only consider
(a): $|\tau-\xi^3| \gtrsim  |\xi|^3$. First we have (since
$|\xi|\leq 2|\xi_1|$),
\begin{equation}\label{4.7}
\displaystyle\int_{D_{42}}
\lesssim
N^{-1}\displaystyle\int_{D_{42}}|\xi|
\widehat{f_\delta}(\xi,\tau)\>\hat{g}(\xi_1,\tau_1)\>\hat{h}(\xi_2,\tau_2).
\end{equation}
Further,
\begin{eqnarray*}
   (\ref{4.7})
 &\lesssim&
   N^{-1}\displaystyle\int_{D_{42}}
   \langle\tau-\xi^3\rangle^{\frac{1}{3}}
   \widehat{f_\delta}(\xi,\tau)\>\hat{g}(\xi_1,\tau_1)\>\hat{h}(\xi_2,\tau_2)\\
 &\lesssim&
   N^{-1}\delta^{a}\>\|f\|_{Y_{0,1-c}}\>\|g\|_{X_{0,b}}\>\|h\|_{X_{0,b}^-},
\end{eqnarray*}
by the same way used in $\displaystyle\int_{D_{12}}$
in the last step.
\hfill$\Box$
\end{proof}

\subsection{Some Variant Local Well-posedness}

Now we turn to obtain a variant local well-posedness result.
Compared with the standard local well-posedness result Theorem 1.1,
it is established in order to fit the $I$-method.
It gives the estimates on the lifetime and the solutions under the
$X_{s,\frac{1}{2}+}^\delta$-norm, with the operator $I_{N,s}$.
Indeed, along the lines of \cite{CKSTT3} and
the estimates from Lemma 3.1--lemma 3.4,
we have the following result as an adaptation of Theorem 1.1
(see \cite{CKSTT2}, \cite{P} cf.).
\begin{cor}
Let $s>0, I=I_{N,s}$, then the solutions obtained in Theorem 1.1 for the initial
data $(u_0,v_0)\in H^s(\R)\times H^s(\R)$ exist on $[-\delta_0,\delta_0]$
with
$$
\delta_0\sim (\|Iu_0\|_{H^1}+\|Iv_0\|_{H^1})^{-\mu};\quad
\|Iu\|_{X_{1,\frac{1}{2}+}^{\delta_0}}+\|Iv\|_{Y_{1,\frac{1}{2}+}^{\delta_0}}\lesssim
\|Iu_0\|_{H^1}+\|Iv_0\|_{H^1},
$$
for some $\mu>0$.
\end{cor}

But we have no intention of exploiting it as our basic of the iteration
to establish the global well-posedness results. In order to extend the lifetime,
we shall reconstruction
it and ultimately establish the refined local result as follows.
\begin{prop}
Let $s>0, I=I_{N,s}$, then
the solutions obtained in Corollary 4.7 exists on $[-\delta,\delta]$
with
\begin{equation}\label{lifetime}
\delta\sim (\|Iu_0\|_{H^1}+\|Iv_0\|_{H^1})^{-2-},
\end{equation}
when $N\gg N_0$ for some large number $N_0$ such that
\begin{equation}\label{N0}
\gamma N_0^{-2}(\|I_{N_0,s}u_0\|_{H^1}+\|I_{N_0,s}v_0\|_{H^1})\sim 1.
\end{equation}
Moreover, the solutions satisfy
\begin{equation}\label{sol}
\|Iu\|_{X_{1,\frac{1}{2}+}^\delta}+\|Iv\|_{Y_{1,\frac{1}{2}+}^\delta}\lesssim
\|Iu_0\|_{H^1}+\|Iv_0\|_{H^1}.
\end{equation}
\end{prop}

In the following text, we may assume that $\|u(t)\|_{L^2_{x}}\sim 1$
by fixing $u_0$ and the $L^2$-mass conservation:
$\|u(t)\|_{L^2_{x}}=\|u_0\|_{L^2_{x}}$. Further, by local result in
Corollary 4.7 and the iteration, one can conclude the existence of
solutions on $[-\delta,\delta]$ in Proposition 4.8 by the estimate
in (\ref{sol}), which implies \emph{a priori} estimate of the
solutions in $H^s(\R)\times H^s(\R)$. Therefore, to prove
Proposition 4.8, we may assume at the beginning that the solutions
exist on the time interval $[-\delta,\delta]$ with the $\delta$
defined in (\ref{lifetime}), and turn to prove (\ref{sol}).
\begin{lem}
Let $s>0$, assume that $(u,v)$ are the solutions of (\ref{NLS-KdV})
on $[-\delta,\delta]$ for small $\delta>0$ with the initial data
$(u_0,v_0)\in H^s(\R)\times H^s(\R)$, then if it satisfies that
\begin{equation}\label{condition}
\alpha \delta^{\frac{13}{16}-}\cdot\|v\|_{Y_{0,\frac{1}{2}+}}^\delta\leq \epsilon_0
\end{equation}
for some small $\epsilon_0$, we have
\begin{equation}\label{X0b}
\|u\|_{X_{0,\frac{1}{2}+}^\delta}\lesssim \|u_0\|_{L^2}.
\end{equation}
\end{lem}
\begin{proof}
Define the operator
\begin{equation}\label{4.9}
\Phi(u)(x,t)=\psi(t)S(t)u_0-i\psi(t)\displaystyle\int_0^t
  S(t-t')\psi(t'/\delta)\left[\alpha (uv)(t')+\beta(|u|^2u)(t')\right]\,dt'.
\end{equation}
Taking $X_{0,\frac{1}{2}+}^{\delta}$ on the two sides of (\ref{4.9}), and by
Lemma 2.4, Lemma 3.1, Lemma 4.3 (when $s=0$), we have
\begin{eqnarray*}
  \|\Phi(u)\|_{X_{0,\frac{1}{2}+}^\delta}
  &\lesssim&
  \|u_0\|_{L^2}+\alpha\|\psi(t/\delta)\>\tilde{u} \tilde{v}\|_{X_{0,-\frac{1}{2}+}}
  + \beta \left\|\psi(t/\delta)\>|\tilde{u}|^2 \tilde{u}
  \right\|_{X_{0,-\frac{1}{2}+}}\\
  &\lesssim&
  \|u_0\|_{L^2}+\alpha\delta^{\frac{13}{16}-}\>\|u\|_{X_{0,\frac{1}{2}+}^\delta}
  \>\|v\|_{Y_{0,\frac{1}{2}+}^\delta}
  + \beta \delta^{\frac{1}{2}-}\>\|u\|_{X_{0,\frac{1}{2}+}^\delta}^3,
\end{eqnarray*}
where $\|u\|_{X_{0,b}^\delta}=\|\tilde{u}\|_{X_{0,b}}$,
$\|v\|_{Y_{0,c}^\delta}=\|\tilde{v}\|_{Y_{0,c}}$.
Thus, by (\ref{condition}), we have
\begin{equation}\label{4.10}
\|\Phi(u)\|_{X_{0,\frac{1}{2}+}^\delta}
\leq
c\|u_0\|_{L^2}+C
\left(\epsilon_0\|u\|_{X_{0,\frac{1}{2}+}^\delta}+
\beta \delta^{\frac{1}{2}-}\>\|u\|_{X_{0,\frac{1}{2}+}^\delta}^3\right)
\end{equation}
for some large constants $c,C>0$. Let the ball
$B\in X_{0,\frac{1}{2}+}^\delta$ be defined as
$$
B=
\left\{u\in X_{0,\frac{1}{2}+}^\delta:
\|u\|_{X_{0,\frac{1}{2}+}^\delta}\leq 2c\|u_0\|_{L^2}\right\},
$$
then by (\ref{4.10}), we have $\Phi$ maps $B$ into itself. We also have
the contraction of $\Phi$ by a similar way.
Thus we complete the proof of the lemma by the fixed point
theory.  \hfill$\Box$
\end{proof}

{\it\noindent Proof of Proposition 4.8. }
By Duhamel's formula and acting the operator $I=I_{N,s}$, we have,
for $t\in [-\varrho,\varrho]$,
\begin{eqnarray*}
  Iu(x,t) &=& \psi(t)S(t)Iu_0-i\psi(t)\displaystyle\int_0^t
  S(t-t')\psi(t'/\varrho)\left[\alpha I(uv)(t')+\beta I(|u|^2u)(t')\right]\,dt', \\
  Iv(x,t) &=& \psi(t)W(t)Iv_0+\psi(t)\displaystyle\int_0^t
  W(t-t')\psi(t'/\varrho)\left[\gamma \partial_xI(|u|^2)(t')
  -\frac{1}{2}\partial_x I(v^2)(t')\right]\,dt'.
\end{eqnarray*}
Therefore, by Lemma 2.4, Corollarys 4.2, 4.4 and Lemmas 4.5, 4.6, we have
\begin{eqnarray}
  \|Iu\|_{X_{1,\frac{1}{2}+}^\varrho}
&\lesssim &
  \|Iu_0\|_{H^1}
  +\alpha \left\|\psi(t/\varrho)I(\tilde{u}\tilde{v})
  \right\|_{X_{1,-\frac{1}{2}+}}
  +\beta \left\|\psi(t/\varrho)I(|\tilde{u}|^2\tilde{u})\right\|_{X_{1,-\frac{1}{2}+}}
  \nonumber\\
&\leq &
  c\|Iu_0\|_{H^1}
  +C\alpha \varrho^{\frac{13}{16}-}\>\|Iu\|_{X_{1,\frac{1}{2}+}^\varrho}
    \>\|Iv\|_{Y_{1,\frac{1}{2}+}^\varrho}
  \nonumber\\
&&
  +C\beta \varrho^{\frac{1}{2}-}\>\|Iu\|_{X_{1,\frac{1}{2}+}^\varrho}
    \>\|u\|_{X_{0,\frac{1}{2}+}^\varrho}^2;\label{4.15}\\
  \|Iv\|_{Y_{1,\frac{1}{2}+}^\varrho}
&\lesssim &
  \|Iv_0\|_{H^1}
  +\gamma\left\|\psi(t/\varrho)\partial_xI\left(|\tilde{u}|^2\right)
  \right\|_{Y_{1,-\frac{1}{2}+}}
  +\left\|\psi(t/\varrho)\partial_xI\left(\tilde{v}^2\right)
  \right\|_{Y_{1,-\frac{1}{2}+}}
  \nonumber\\
&\leq &
  c\|Iv_0\|_{H^1}
  +C\gamma\Big(\varrho^{a}\>\|Iu\|_{X_{1,\frac{1}{2}+}^\varrho}
    \>\|u\|_{X_{0,\frac{1}{2}+}^\varrho}
    \nonumber\\
&&
  +(N^{-1}\varrho^{a}+\varrho^{\frac{13}{16}-})
    \>\|Iu\|_{X_{1,\frac{1}{2}+}^\varrho}^2\Big)
    +C\varrho^{\frac{1}{2}-}\>\|Iv\|_{Y_{1,\frac{1}{2}+}^\varrho}^2
    \label{4.16}
\end{eqnarray}
for some constants $c,C>0$. Let $\delta$ be the quantity satisfying
\begin{equation}\label{4.17}
\begin{array}{c}
\alpha \delta^{\frac{13}{16}-}R(N)\leq \epsilon_0;\quad
\beta \delta^{\frac{1}{2}-}\>\|u_0\|_{L^2}^2\leq \epsilon_0;\quad\\
\gamma\delta^{a}\>\|u_0\|_{L^2}\leq \epsilon_0;\quad
\gamma N^{-1}\delta^{a}\>R(N)\leq \epsilon_0;\quad
\gamma\delta^{\frac{13}{16}-}R(N)\leq \epsilon_0;\quad
\delta^{\frac{1}{2}-}R(N)\leq \epsilon_0
\end{array}
\end{equation}
for some small $\epsilon_0$ and $R(N)=\|Iu_0\|_{H^1}+\|Iv_0\|_{H^1}$.
We claim that for any $\varrho\in[0,\delta]$,
\begin{equation}\label{claim}
\|Iu\|_{X_{1,\frac{1}{2}+}^\varrho}+\|Iv\|_{Y_{1,\frac{1}{2}+}^\varrho}\leq
2c R(N).
\end{equation}
Indeed, it can be shown by the iteration which we present as follows.
We only consider the positive time,
since it is similar to the negative time.
First, by Corollary 4.7, we have
$$
\|Iu\|_{X_{1,\frac{1}{2}+}^{\delta_0}}
+\|Iv\|_{Y_{1,\frac{1}{2}+}^{\delta_0}}\leq
2c R(N),
$$
then by Lemma 2.4 (i), we obtain
$$
R_1(N)\equiv\|Iu(\delta_0)\|_{H^1}+\|Iv(\delta_0)\|_{H^1}\leq 2c\bar{c}R(N)
$$
for some constant $\bar{c}>0$.
Now we take $(u(\delta_0),v(\delta_0))$ for the new initial data, and
employ Corollary 4.7 again, then we obtain that
for some $\delta_1\sim \delta_0$,
$$
\|Iu\|_{X_{1,\frac{1}{2}+}^{[\delta_0,\delta_1]}}
+\|Iv\|_{Y_{1,\frac{1}{2}+}^{[\delta_0,\delta_1]}}
\leq
2c R_1(N)
\leq
4c^2\bar{c}R(N).
$$
Since
$\|f\|_{X_{s,b}^{[0,\delta_1]}(\phi)}\leq
\|f\|_{X_{s,b}^{[0,\delta_0]}(\phi)}+
\|f\|_{X_{s,b}^{[\delta_0,\delta_1]}(\phi)}$ by Lemma 6.2, we get
\begin{equation}\label{419}
\|Iu\|_{X_{1,\frac{1}{2}+}^{\delta_1}}
+\|Iv\|_{Y_{1,\frac{1}{2}+}^{\delta_1}}
\leq
2c (1+2c\bar{c})R(N).
\end{equation}
By (\ref{4.17}) and (\ref{419}), we have (\ref{condition}), and thus we
obtain (\ref{X0b}) by Lemma 4.9. Therefore, inserting (\ref{X0b}) and
(\ref{419}) into  (\ref{4.15}) and  (\ref{4.16}), we have exactly
$$
\|Iu\|_{X_{1,\frac{1}{2}+}^{\delta_1}}
+\|Iv\|_{Y_{1,\frac{1}{2}+}^{\delta_1}}
\leq
2c R(N).
$$
The process above can always be repeated under (\ref{4.17}),
and ultimately we prove the claim (\ref{claim}) and obtain the proposition.
\hfill$\Box$

\subsection{The Global Well-posedness}
In this section, we establish the global result by combining
Proposition 4.8 and the results in \cite{P}.
We only consider the positive time in the following.
Compared to
the precess in Section 6 in \cite{P},
it just needs to modify the estimate on the lifetime.
Define
\begin{eqnarray*}
  M_I(t)
 &=& \|Iu\|_{L^2}; \\
  L_I(t)
&=&
  \alpha\|Iv\|_{L^2}^2+2\gamma \displaystyle\int\! \mbox{Im}(Iu
  \>I\bar{u}_x)\,dx; \\
  E_I(t)
&=&
  \alpha\gamma \displaystyle\int Iv
  \>|Iu|^2\,dx+\gamma\|Iu_x\|_{L^2}^2+\frac{\alpha}{2}\|Iv_x\|_{L^2}^2
  -\frac{\alpha}{6}\|Iv\|_{L^3}^3+\frac{\beta\gamma}{2}\|Iv\|_{L^4}^4,
\end{eqnarray*}
then by the Sobolev interpolation inequalities, we have (see \cite{P}
for the details)
\begin{equation}\label{4.19}
\|Iu\|_{H^1}^2+\|Iv\|_{H^1}^2
\lesssim
|E_I|+|L_I|^{\frac{5}{3}}+|M_I|^8+1.
\end{equation}
Moreover, we have the following estimates, which are proved in \cite{P}.
\begin{lem}
Let $I=I_{N,s}, s>\dfrac{1}{2}$, $(u,v)$ is the solution of (\ref{NLS-KdV}),
then
\begin{eqnarray*}
|E_I(\delta)-E_I(0)|
\lesssim
\left(N^{-1+}\delta^{\frac{1}{2}-}+N^{\frac{7}{4}+}\right)
\left(\|Iu\|_{X_{1,\frac{1}{2}+}}^3
+\|Iv\|_{Y_{1,\frac{1}{2}+}}^3\right)+N^{-2+}\cdot\qquad\\
\left(\|Iu\|_{X_{1,\frac{1}{2}+}}^4+\|Iv\|_{Y_{1,\frac{1}{2}+}}^4\right)
+N^{-3+}\|Iu\|_{X_{1,\frac{1}{2}+}}^4
\left(\|Iu\|_{X_{1,\frac{1}{2}+}}^2+\|Iv\|_{Y_{1,\frac{1}{2}+}}\right).
\end{eqnarray*}
\end{lem}
\begin{lem}
Let $I, s, (u,v)$ be the same with Lemma 4.10,
then
$$
|L_I(\delta)-L_I(0)|
\lesssim
N^{-2+}\delta^{\frac{1}{2}-}
\left(\|Iu\|_{X_{1,\frac{1}{2}+}}^3+\|Iv\|_{Y_{1,\frac{1}{2}+}}^3\right)
+N^{-3+}\|Iu\|_{X_{1,\frac{1}{2}+}}^4.
$$
\end{lem}
Further, we have the trivial estimate of $M_I(t)$ that
$$
M_I(t)\lesssim \|u(t)\|_{L^2}\sim 1,
$$
which follows from
the $L^2$-mass conservation and $m(\xi)\leq1$.

Fix the large number $N$, $s>\dfrac{1}{2}$. By Proposition 4.8,
the solution $(u,v)$ of the system (\ref{NLS-KdV}) exists on $[0,\delta]$,
with
$$
\delta\sim (\|Iu_0\|_{H^1}+\|Iv_0\|_{H^1})^{-2-}\gtrsim N^{-2(1-s)-}
$$
by (\ref{II}).
We repeat this local existence results by iteration. In order to ensure
the same length of the lifespan, we need to get the uniform control
of $H^1$-norm of the solution at $t=k\delta$ for $k=1,2,\cdots$, which
follows from the uniform control of $|E_I|$ and $|L_I|$. More precisely,
we shall obtain that
\begin{equation}\label{4.20}
|E_I(k\delta)|\leq cN^{2(1-s)};\quad
|L_I(k\delta)|\leq cN^{1-s},
\end{equation}
which imply by (\ref{4.19}) that
$\|Iu(k\delta)\|_{H^1}^2+\|Iv(k\delta)\|_{H^1}^2\leq \bar{c}N^{2(1-s)}$
for the constants $c,\bar{c}$ independent of $k,N$. We note that the
condition (\ref{N0}) is valid in every step. By Lemmas 4.10, 4.11, and
the estimate of (\ref{sol}) in each step, we have
$$
\begin{array}{c}
|E_I(k\delta)-E_I(0)|
\leq \tilde{c}k\left(
\left(N^{-1+}\delta^{\frac{1}{2}-}+N^{-\frac{7}{4}+}\right)N^{3(1-s)}
+N^{-2+}N^{4(1-s)}+N^{-3+}N^{6(1-s)}\right);\\
|L_I(k\delta)-L_I(0)|
\leq \tilde{c}k\left(N^{-2+}\delta^{\frac{1}{2}-}N^{3(1-s)}
+N^{-3+}N^{4(1-s)}\right)
\end{array}
$$
for the constant $\tilde{c}$ independent of $k,N$.

Set $T=k\delta$, for (\ref{4.20}), we only need to show
\begin{eqnarray}
T\delta^{-1}\left(\left(N^{-1+}\delta^{\frac{1}{2}-}
+N^{-\frac{7}{4}+}\right)N^{3(1-s)}
+N^{-2+}N^{4(1-s)}+N^{-3+}N^{6(1-s)}\right)
\lesssim
N^{2(1-s)};\label{4.21}\\
T\delta^{-1}\left(N^{-2+}\delta^{\frac{1}{2}-}N^{3(1-s)}
+N^{-3+}N^{4(1-s)}\right)
\lesssim
N^{1-s}.\qquad\qquad\qquad\label{4.22}
\end{eqnarray}
In order to $T\sim N^{0+}$, and note that $\delta\sim N^{-2(1-s)-}$,
(\ref{4.21}) is fulfilled if
$$
N^{-1+}N^{-(1-s)+}N^{3(1-s)}
+N^{-\frac{7}{4}+}N^{3(1-s)}
+N^{-2+}N^{4(1-s)}+N^{-3+}N^{6(1-s)}
\lesssim
N^{0-},
$$
which is valid if
$$
-1-(1-s)+3(1-s)<0;\,
-\dfrac{7}{4}+3(1-s)<0;\,
-2+4(1-s)<0;\,
-3+6(1-s)<0,
$$
which hold when $s>\dfrac{1}{2}$. Similarly, (\ref{4.22}) is fulfilled if
$$
N^{-2+}N^{-(1-s)+}N^{3(1-s)}
+N^{-3+}N^{4(1-s)}
\lesssim
N^{(s-1)-}.
$$
It is valid if
$$
-2-(1-s)+3(1-s)<s-1;\,
-3+4(1-s)<s-1,
$$
they hold when $s>\dfrac{2}{5}$.
Therefore, we prove the global well-posedness in $H^s(\R)\times H^s(\R)$
when $s>\dfrac{1}{2}$ and thus finish the proof of Theorem 1.3.

 \vspace{0.3cm}
\section{The Proof of Theorem 1.2}

Suppose for the contradiction that the system (\ref{NLS-KdV}) is
locally well-posed on $[0,\delta]$ for $\delta\in (0,1)$,
and the solution map $(u_0,v_0)\mapsto (u,v)$ is $C^2$ from
$H^{s}(\R)\times H^{l}(\R)$ to
$C_t^0([0,\delta];H^{s}(\R)\times H^{l}(\R))$.
Then, by the Picard iterative
scheme, so is the operator
$A=(A_1,A_2): H^{s}(\R)\times H^{l}(\R)\rightarrow
C_t^0([0,\delta];H^{s}(\R)\times H^{l}(\R))$ defined as
\begin{eqnarray*}
A_1(u_0,v_0)&=&-i\displaystyle\int_0^t
  S(t-t')\left[\alpha (S(t')u_0\cdot W(t')v_0)
  +\beta\left(|S(t')u_0|^2S(t')u_0\right)\right]\,dt';\\
A_2(u_0,v_0)&=& \displaystyle\int_0^t
  W(t-t')\left[\gamma \partial_x(|S(t')u_0|^2)(t')
  -\frac{1}{2}\partial_x(W(t')v_0)^2\right]\,dt'.
\end{eqnarray*}
In particular, $A_2$ is $C^2$-differentiable from
$H^{s}(\R)\times H^{l}(\R)$ to
$C_t^0([0,\delta];H^{l}(\R))$.

Fix a large number $N\gg1$ such that $(N-\dfrac{1}{2})^3=2k\pi$
for some $k\in \mathbb{N}$,
and let the sets
\begin{eqnarray*}
  \Upsilon   &=& \left\{\xi\in \R:
     \left|(\xi+\frac{1}{2})-N\right|\leq \frac{1}{100N^2}\right\}; \\
  \Upsilon_1 &=& \left\{\xi\in \R:
     \left|\xi-\left(N-\frac{1}{2}N^2
     -\frac{3}{8}\right)\right|\leq \frac{1}{N}\right\}; \\
  \Upsilon_2 &=& \left\{\xi\in \R:
     \left|(2\xi+\frac{1}{4})-N^2\right|\leq \frac{1}{100N}\right\}; \\
  \Lambda   &=& \left\{\xi\in \R: |\xi-1|\leq \frac{1}{N^n}\right\}
\end{eqnarray*}
for some $n\in \mathbb{N}$.
Note that
\begin{equation}\label{5.1}
\left\{(\xi_1,\xi_2)\in \R^2:
\xi=\xi_1+\xi_2, \xi\in \Upsilon, \xi_2\in \Upsilon_2\right\}
\subset
\left\{(\xi_1,\xi_2)\in \R^2: \xi=\xi_1+\xi_2,\xi_1\in \Upsilon_1\right\}.
\end{equation}
Put the initial data $(u_0,v_0)$ such that
$$
\widehat{u_0}(\xi)= \epsilon_0N^{-2s+\frac{1}{2}} \chi_{\Upsilon_1}(\xi);\quad
\widehat{v_0}(\xi)= \epsilon_0N^{\frac{n}{2}} \chi_{\Lambda}(\xi),
$$
then $\|u_0\|_{H^{s}},\|v_0\|_{H^{l}}\sim \epsilon_0$. We may set $\delta=1$
by choosing
$\epsilon_0$ small enough.

Further,
$\|A_2\|_{C_t^0([0,1];H^{s})}$ is equal to
\begin{eqnarray}
  &&
   \sup_{0\leq t \leq 1}\left\|\displaystyle\int_0^t
   W(t-t')\left[\gamma \partial_x(|S(t')u_0|^2)
   -\frac{1}{2}\partial_x(W(t')v_0)^2\right]\,dt'\right\|_{H^{l}}\nonumber\\
  &=&
   \sup_{0\leq t \leq 1}\biggl\|\xi\langle\xi\rangle^{l}
   \displaystyle\int_0^t\!\!\!\int
   \exp\left\{i(t-t')\xi^3\right\}
   \big[\gamma \exp\left\{-it'(\xi-\xi_2)^2\right\}
   \exp\left\{it'\xi_2^2\right\}
   \widehat{u_0}(\xi-\xi_2)\widehat{u_0}(\xi_2)\nonumber\\
  &&
   \qquad-\frac{1}{2}
   \exp\left\{it'(\xi-\xi_2)^3\right\}\exp\left\{it'\xi_2^3\right\}
   \widehat{v_0}(\xi-\xi_2)\widehat{v_0}(\xi_2)\big] \,d\xi_2 dt'
   \biggl\|_{L^2_\xi}\nonumber\\
  &\gtrsim&
   \gamma\biggl\|\xi\langle\xi\rangle^{l}\!\!
   \displaystyle\int_0^{1}\!\!\!\!\int
   \exp\left\{i\xi^3\right\}
   \exp\left\{-it'\xi(\xi^2+\xi-2\xi_2)\right\}
   \widehat{u_0}(\xi-\xi_2)\widehat{u_0}(\xi_2)
   \,d\xi_2 dt'\biggl\|_{L^2_\xi}\nonumber\\
  &&
   \qquad-\biggl\|\xi\langle\xi\rangle^{l}\!\!
   \displaystyle\int_0^{1}\!\!\!\!\int
   \exp\left\{i\xi^3\right\}\exp\left\{-3it'\xi(\xi-\xi_2)\xi_2\right\}
   \widehat{v_0}(\xi-\xi_2)\widehat{v_0}(\xi_2)
   \,d\xi_2 dt'\biggl\|_{L^2_\xi}.\label{5.2}
\end{eqnarray}
The first term of (\ref{5.2}) has a lower bound of
\begin{equation}\label{5.3}
    \gamma N^{l+1}\biggl\|
   \displaystyle\int_0^{1}\!\!\!\!\int
   \exp\left\{i\xi^3\right\}
   \exp\left\{-it'\xi(\xi^2+\xi-2\xi_2)\right\}
   \widehat{u_0}(\xi-\xi_2)\widehat{u_0}(\xi_2)
   \,d\xi_2 dt'\biggl\|_{L^2_\xi(\Upsilon)}.
\end{equation}
When $\xi\in \Upsilon$, set $\xi=N-\dfrac{1}{2}+p(\xi,N)$,
then $|p|\leq \dfrac{1}{100N^2}$. Therefore,
$$
\left|\xi^3-\left(N-\frac{1}{2}\right)^3\right|
\leq 2 N^2\cdot|p|\leq \dfrac{1}{50}.
$$
Besides, when $\xi\in \Upsilon,\xi_2\in \Upsilon_2$, we have
$$
\left|\xi(\xi^2+\xi-2\xi_2)\right|\leq\dfrac{1}{10}.
$$
Note that $(N-\dfrac{1}{2})^3=2k\pi$,
therefore, for any $0\leq t'\leq 1$,
we have,
\begin{equation}\label{5.4}
\mbox{Re}\left(\exp\left\{i\xi^3\right\}
   \exp\left\{-it'\xi(\xi^2+\xi-2\xi_2)\right\}\right)>\frac{1}{2}.
\end{equation}
Thus, by (\ref{5.1}) and (\ref{5.4}), we obtain
\begin{eqnarray*}
(\ref{5.3})
  &\gtrsim &
   \gamma \epsilon_0^2 N^{l+1}N^{-4s+1}\>m(\Upsilon_2)\>m(\Upsilon)^{\frac{1}{2}}
  \sim
   \gamma \epsilon_0^2 N^{l-4s},
\end{eqnarray*}
where $m(\cdot)$ is the Lebesgue measure.

Second term of (\ref{5.2}), by the support of $v_0$, have a
bound of
$$
\left\| \displaystyle\int_0^{1}\!\!\!\!\int
   \widehat{v_0}(\xi-\xi_2)\widehat{v_0}(\xi_2)
   \,d\xi_2 dt'\right\|_{L^2_\xi(\{|\xi-2|\leq \frac{2}{N^{n}}\})}
\lesssim
\epsilon_0^2 N^{-\frac{n}{2}}.
$$
Therefore, by choosing $n$ large enough, we have
\begin{equation}\label{55}
\|A_2(u_0,v_0)\|_{C_t^0([0,1];H^{s})}
\gtrsim
\gamma \epsilon_0^2 N^{l-4s}.
\end{equation}
Since $A_2$ is $C^2$-differentiable, we must have
$$
\|A_2(u_0,v_0)\|_{C_t^0([0,1];H^{s})}
\lesssim \|u_0\|_{H^{s}}^2+\|v_0\|_{H^{l}}^2,
$$
but it fails to hold when $l>4s$ by (\ref{55}).
This completes the proof of Theorem 1.2.

{\noindent\it Remark.} More facts related to the condition $l<4s$
may be interesting to the readers. As we see, the condition appears in the
bilinear estimate in Lemma 3.4,
which is necessary in the framework of Bourgain method.
But on the other hand, it is optimal. More precisely, if
$l>4s$, then for any $b_1,b_2,b_3\in\R$, the estimates
$$
    \|\partial_x(u_1\overline{u_2})\|_{Y_{l,b_1}}
    \lesssim
    \|u_1\|_{X_{s,b_2}}\>\|u_2\|_{X_{s,b_3}}
$$
fail to hold. The proof is based on the counterexample of $u_1,u_2$ such that
$$
\widehat{u_1}(\xi,\tau)
=\chi_{\Upsilon_1}\cdot\chi_{\{|\tau-\xi^2|\leq 100\}}(\xi,\tau);\quad
\widehat{u_2}(\xi,\tau)
=\chi_{\Upsilon_2}\cdot\chi_{\{|\tau-\xi^2|\leq 10\}}(\xi,\tau),
$$
and specially take the integration in the left hand side over
$$
\{(\xi,\tau): \xi\in \Upsilon, |\tau-\xi^2|\leq 10\}
$$
in the spacetime-frequency space. The detailed computation is omitted here.

 \vspace{0.3cm}

\section{Appendix: Some Auxiliary Lemmas}

As some handy tools, we give some properties and estimates on the restricted spaces
$X_{s,b}(\phi)$ in the following.
\begin{lem}
For the time interval $\Omega$ and the function $f\in X_{s,b}^{\Omega}(\phi)$,
there exists an extension
$\tilde{f}\in X_{s,b}(\phi)$ such that $\tilde{f}=f$ on
$\Omega$ and
$$
\|f\|_{X_{s,b}^{\Omega}(\phi)}
=\left\|\tilde{f}\right\|_{X_{s,b}(\phi)}.
$$
Moreover, it holds that for any $s'\leq s$,
\begin{equation}\label{L4.1}
\|I_{N,s}f\|_{X_{1,b}^{\Omega}(\phi)}
=\left\|I_{N,s}\tilde{f}\right\|_{X_{1,b}(\phi)};
\quad
\|f\|_{X_{s',b}^{\Omega}(\phi)}
=\left\|\tilde{f}\right\|_{X_{s',b}(\phi)}.
\end{equation}
\end{lem}
\begin{proof}
Fix the function $f\in X_{s,b}^{\Omega}(\phi)$, and set
$$
M_{s,b}=\left\{F\in X_{s,b}(\phi):
F\big|_{t\in\Omega}=f\big|_{t\in\Omega}\right\},
$$
then $ \|f\|_{X_{s,b}^{\Omega}(\phi)}=\inf\limits_{F\in
M_{s,b}}\|F\|_{X_{s,b}(\phi)}.$ Since $X_{s,b}$ is a Hilbert space
and $M_{s,b}$ is a closed convex subset of $X_{s,b}$,
there exists a minimum $\tilde{f}$ in $M_{s,b}$. Move precisely, note that
$$
M_{s,b}= \overset{\circ}{M}_{s,b}+\{F\}
\equiv \{g+F:g\in \overset{\circ}{M}_{s,b}\}
$$
for any $F\in M_{s,b}$, where $\overset{\circ}{M}_{s,b}$ is the closed linear
subspace of $X_{s,b}$ defined as
$$
\overset{\circ}{M}_{s,b}=\{g\in X_{s,b}(\phi):
g|_{\Omega}=0\},
$$
therefore,
$\tilde{f}$ is exactly the function that
$$
\tilde{f}\bot \overset{\circ}{M}_{s,b} \mbox{\quad in\quad}  X_{s,b}(\phi).
$$
Therefore, for (\ref{L4.1}), we only need to show that
$$
I_{N,s}\tilde{f}\bot \overset{\circ}{M}_{1,b} \mbox{\quad in\quad}  X_{1,b}(\phi);\quad
\tilde{f}\bot \overset{\circ}{M}_{s',b} \mbox{\quad in\quad}  X_{s',b}(\phi)
$$
for $s'\leq s$, but it is obvious.  \hfill$\Box$
\end{proof}

Now we show that the spaces $X_{s,b}(\phi)$ have the property of
the norm-subadditivity about the restricted domain.
\begin{lem}
For every $f\in X_{s,b}^{\Omega}(\phi)$,
it holds for any $\Omega_0\subset \Omega$ that
$$
\|f\|_{X_{s,b}^{\Omega}(\phi)}\leq
\|f\|_{X_{s,b}^{\Omega_0}(\phi)}+
\|f\|_{X_{s,b}^{\Omega/\Omega_0}(\phi)}.
$$
\end{lem}
\begin{proof}
By Lemma 6.1, there exists a
function $\tilde{f}\in X_{s,b}(\phi)$ such that $\tilde{f}=f$ on
$\Omega$ and
$$
\|f\|_{X_{s,b}^{\Omega}(\phi)}
=\left\|\tilde{f}\right\|_{X_{s,b}(\phi)}.
$$
Define the operator $P_b$ and its inverse operator $P_b^{-1}$ as
$$
\widehat{P_b f}(\xi,\tau)=\langle\tau+\phi(\xi)\rangle^{2b}\hat{f}(\xi,\tau),
\quad
\widehat{P_b^{-1} f}(\xi,\tau)=\langle\tau+\phi(\xi)\rangle^{-2b}\hat{f}(\xi,\tau).
$$
We claim that
$$
\left\|\tilde{f}\right\|_{X_{s,b}(\phi)}
=\left\|P_b^{-1} \left(P_b\tilde{f}\cdot j\right)\right\|_{X_{s,b}(\phi)}
=\left\|P_b^{-\frac{1}{2}} \left(P_b\tilde{f}\cdot j\right)\right\|_{L^2_t H^s_x}
$$
for any $j(t)\in C_c^\infty(\R)$, such that
$$
j(t)\equiv 1 \mbox{\,\, on \,\,} \Omega'\supset\supset \Omega.
$$
Indeed, by the proof of Lemma 6.1, one only needs to show that
$$
P_b^{-1} \left(P_b\tilde{f}\cdot j\right)\in M_{s,b};\quad
P_b^{-1} \left(P_b\tilde{f}\cdot j\right)
\bot \overset{\circ}{M}_{s,b} \mbox{\quad in\quad}  X_{s,b}(\phi).
$$
The first term is easy to check and we omit the details.
On the other hand, for any $g\in \overset{\circ}{M}_{s,b}$, we have
\begin{eqnarray*}
  \left\langle P_b^{-1} \left(P_b\tilde{f}\cdot j\right), g\right\rangle
 &=&
  \displaystyle\int\langle\xi\rangle^{2s} \langle\tau+\phi(\xi)\rangle^{2b}
  \F\left(P_b^{-1} \left(P_b\tilde{f}\cdot j\right)\right)(\xi,\tau)
  \cdot\overline{\F g}(\xi,\tau)
  d\xi d\tau\\
 &=&
  \displaystyle\int\langle\xi\rangle^{2s}
  \left(\F\left(P_b\tilde{f}\right)\ast\F j\right)(\xi,\tau)
  \cdot\overline{\F g}(\xi,\tau)
  d\xi d\tau\\
 &=&
  \displaystyle\int\langle\xi\rangle^{2s}
  \F\left(P_b\tilde{f}\right)(\xi,\tau)
  \cdot\overline{\F (gj(-\cdot))}(\xi,\tau)
  d\xi d\tau\\
 &=&
  \left\langle \tilde{f}, gj(-\cdot)\right\rangle
  \\
 &=&0,
\end{eqnarray*}
since $gj(-\cdot)\in \overset{\circ}{M}_{s,b}$,
and $\tilde{f}\bot \overset{\circ}{M}_{s,b}$ in $X_{s,b}(\phi)$,
where $\langle \cdot, \cdot\rangle$ is the inner product in $X_{s,b}(\phi)$.
Therefore, we have the claim.

Furthermore, we keep in mind that
\begin{equation}\label{c}
P_b^{-\frac{1}{2}} \left(P_b\tilde{f}\cdot j\right)\Big|_{t\in\Omega}
=P_b^{-\frac{1}{2}} f\Big|_{t\in\Omega};\quad
P_b^{-\frac{1}{2}} \left(P_b\tilde{f}\cdot j\right)\Big|_{t\in(\mbox{supp}j)^c}
\equiv0.
\end{equation}
For the simplicity, we only prove the special result that
$$
\|f\|_{X_{s,b}^{[0,\delta]}(\phi)}\leq
\|f\|_{X_{s,b}^{[0,\delta_0]}(\phi)}+
\|f\|_{X_{s,b}^{[\delta_0,\delta]}(\phi)}
$$
for $0\leq \delta_0\leq \delta$, which follows from
\begin{equation}\label{5.6}
\|f\|_{X_{s,b}^{[0,\delta]}(\phi)}\leq
\liminf\limits_{\epsilon\rightarrow0}
\|f\|_{X_{s,b}^{[0,\delta]/E_\epsilon}(\phi)}\leq
\|f\|_{X_{s,b}^{[0,\delta_0]}(\phi)}+
\|f\|_{X_{s,b}^{[\delta_0,\delta]}(\phi)},
\end{equation}
where $E_\epsilon=(\delta_0-4\epsilon,\delta_0+4\epsilon)$.

We first show the second inequality of (\ref{5.6}). Indeed, for
$f_1,f_2$ such that $\|f\|_{X_{s,b}^{[0,\delta_0]}(\phi)}
=\left\|\widetilde{f_1}\right\|_{X_{s,b}(\phi)},
\|f\|_{X_{s,b}^{[\delta_0,\delta]}(\phi)}
=\left\|\widetilde{f_2}\right\|_{X_{s,b}(\phi)}$ and the function
$j_1(t),j_2(t)$ such that
$$
\begin{array}{c}
\mbox{supp}j_1\subset (-2\epsilon, \delta_0+2\epsilon),\,\,
j_1(t)\equiv 1 \mbox{\,\,on\,\,}(-\epsilon, \delta_0+\epsilon);\\
\mbox{supp}j_2\subset (\delta_0-2\epsilon, \delta+2\epsilon),\,\,
j_1(t)\equiv 1 \mbox{\,\,on\,\,}(\delta_0-\epsilon, \delta+\epsilon),
\end{array}
$$
we have for any small $\epsilon>0$,
$$
P_b^{-1} \left(P_b\tilde{f_1}\cdot j_1+
P_b\tilde{f_2}\cdot j_2\right)\Big|_{t\in [0,\delta]/E_\epsilon}
=f\big|_{t\in [0,\delta]/E_\epsilon}.
$$
Therefore,
$$
\|f\|_{X_{s,b}^{[0,\delta]/E_\epsilon}(\phi)}\leq
\left\|P_b^{-1} \left(P_b\widetilde{f_1}\cdot j_1+
P_b\widetilde{f_2}\cdot j_2\right)\right\|_{X_{s,b}(\phi)}\leq
\|f\|_{X_{s,b}^{[0,\delta_0]}(\phi)}+
\|f\|_{X_{s,b}^{[\delta_0,\delta]}(\phi)},
$$
since
$$
\left\|P_b^{-1} \left(P_b\widetilde{f_1}\cdot
j_1\right)\right\|_{X_{s,b}(\phi)}
=\|f\|_{X_{s,b}^{[0,\delta_0]}(\phi)},\quad
\left\|P_b^{-1} \left(P_b\widetilde{f_2}\cdot
j_2\right)\right\|_{X_{s,b}(\phi)}
=\|f\|_{X_{s,b}^{[\delta_0,\delta]}(\phi)}.
$$

Now we turn to the first term of (\ref{5.6}). We choose the function
$j(t)$ such that
$$
\mbox{supp}j\subset (-2\epsilon, \delta+2\epsilon),\,\,
j_1(t)\equiv 1 \mbox{\,\,on\,\,}(-\epsilon, \delta+\epsilon),
$$
then for $f,f_\epsilon$ such that
$\|f\|_{X_{s,b}^{[0,\delta]}(\phi)}
=\left\|\tilde{f}\right\|_{X_{s,b}(\phi)},
\|f\|_{X_{s,b}^{[0,\delta]/E_\epsilon}(\phi)}
=\left\|\widetilde{f_\epsilon}\right\|_{X_{s,b}(\phi)}$,
we have
$$
\left\|P_b^{-\frac{1}{2}} \left(P_b\tilde{f}\cdot j\right)\right\|_{L^2_t H^s_x}
=\|f\|_{X_{s,b}^{[0,\delta]}(\phi)};\quad
\left\|P_b^{-\frac{1}{2}} \left(P_b\widetilde{f_\epsilon}\cdot j\right)
\right\|_{L^2_t H^s_x}
=\|f\|_{X_{s,b}^{[0,\delta]/E_\epsilon}(\phi)}.
$$
Moreover, by (\ref{c}) we have,
$$
m\left\{t\in \R:
\left\|P_b^{-\frac{1}{2}} \left(P_b\left(\tilde{f}-\widetilde{f_\epsilon}
  \right)\cdot j\right)(t)\right\|_{H^s_x}>\varepsilon\right\}
\lesssim
\epsilon
$$
for any $\varepsilon>0$, where $m$ is the Lebesgue measure, which implies
that
$$
\left\|P_b^{-\frac{1}{2}} \left(P_b\widetilde{f_\epsilon}\cdot j\right)(t)
\right\|_{H^s_x}
\rightarrow
\left\|P_b^{-\frac{1}{2}} \left(P_b\tilde{f}\cdot j\right)(t)\right\|_{H^s_x}
\mbox{\,\, in measure}.
$$
Thus by Fatou's lemma, we have
$$
\left\|P_b^{-\frac{1}{2}} \left(P_b\tilde{f}\cdot j\right)\right\|_{L^2_t H^s_x}
\leq
\liminf_{\epsilon\rightarrow0}
\left\|P_b^{-\frac{1}{2}} \left(P_b\widetilde{f_\epsilon}\cdot j\right)
\right\|_{L^2_t H^s_x},
$$
and thus complete the proof of the lemma.   \hfill$\Box$
\end{proof}

{\noindent \it Remark.} By the similar argument as above, we can prove that
$
y(\delta)\equiv\|f\|_{X_{s,b}^\delta(\phi)}
$
is left continuous for every fixed function $f$.

\begin{lem}Let $\delta\in (0,1)$, $b,b'\in [0,\dfrac{1}{2})$ with $b'\geq b$,
then
$$
 \|f\|_{X_{s,b}^\delta(\phi)}\lesssim \delta^{b'-b}
 \|f\|_{X_{s,b'}^\delta(\phi)}.
$$
\end{lem}
\begin{proof}
Let $\tilde{f}$ be the extension of
$f\in X_{s,b'}^\delta(\phi)$ such that $\tilde{f}=f$ on
$[-\delta,\delta]$ and
$
\|f\|_{X_{s,b'}^{\delta}(\phi)}
=\left\|\tilde{f}\right\|_{X_{s,b'}(\phi)}.
$
By Lemma 2.4 (iv),
$$
\|f\|_{X_{s,b}^\delta(\phi)}\leq
\left\|\psi(t/\delta)\tilde{f}\right\|_{X_{s,b}(\phi)}\lesssim\delta^{b'-b}
\left\|\tilde{f}\right\|_{X_{s,b'}(\phi)}
=\delta^{b'-b}\|f\|_{X_{s,b'}^\delta(\phi)}.
$$
This completes the proof of the lemma.  \hfill$\Box$
\end{proof}

\vspace{0.3cm}
{\noindent\bf Acknowledgment:} The author would like to express the appreciation
to his advisor Professor Yongsheng Li for his helpful conservations and encouragement.

\end{document}